\documentclass[11pt,a4paper]{amsart}
\usepackage[english]{babel}
\usepackage{graphicx}
\usepackage{framed}
\usepackage[normalem]{ulem}
\usepackage{amsmath}
\usepackage{amsthm}
\usepackage{amssymb}
\usepackage{amsfonts}
\usepackage[foot]{amsaddr}
\usepackage{mathrsfs,amsmath}
\usepackage{enumerate}
\usepackage[utf8]{inputenc}
\usepackage[top=1 in,bottom=1in, left=1 in, right=1 in]{geometry}
\usepackage{mathtools}
\usepackage{tikz-cd}
\usepackage[font=small,labelfont=bf]{caption} 
\usepackage[subrefformat=parens]{subcaption}
\usepackage{cite}
\usepackage{abstract}

\usetikzlibrary{
	knots,
	hobby,
	decorations.pathreplacing,
	shapes.geometric,
	calc
}

\tikzset{
	marrow/.style={decoration={markings,mark=at position 0.5 with {\arrow{#1}}}, postaction=decorate}
}

\tikzset{
	knot diagram/every strand/.append style={
		ultra thick,
		black
	},
	show curve controls/.style={
		postaction=decorate,
		decoration={show path construction,
			curveto code={
				\draw [blue, dashed]
				(\tikzinputsegmentfirst) -- (\tikzinputsegmentsupporta)
				node [at end, draw, solid, red, inner sep=2pt]{};
				\draw [blue, dashed]
				(\tikzinputsegmentsupportb) -- (\tikzinputsegmentlast)
				node [at start, draw, solid, red, inner sep=2pt]{}
				node [at end, fill, blue, ellipse, inner sep=2pt]{}
				;
			}
		}
	},
	show curve endpoints/.style={
		postaction=decorate,
		decoration={show path construction,
			curveto code={
				\node [fill, blue, ellipse, inner sep=2pt] at (\tikzinputsegmentlast) {}
				;
			}
		}
	}
}

\usepackage{enumitem}
\setlist{  
	listparindent=\parindent,
	parsep=0pt,
}

\usepackage{bbm}

\usepackage{stackrel}

\newcommand{\Hom}{{\rm Hom}}

\newcommand{\GW}{{\rm GW}}
\newcommand{\cp}{{\rm cp}}

\newcommand{\R}{\mathbb{R}}
\newcommand{\N}{\mathbb{N}}
\newcommand{\Q}{\mathbb{Q}}
\newcommand{\Z}{\mathbb{Z}}
\newcommand{\C}{\mathbb{C}}

\newcommand{\eps}{\varepsilon}

\newtheorem{thm}{Theorem}[section]
\newtheorem{cor}[thm]{Corollary}
\newtheorem{lem}[thm]{Lemma}

\newtheorem{prop}[thm]{Proposition}

\theoremstyle{definition}
\newtheorem{defn}[thm]{Definition}

\theoremstyle{remark}
\newtheorem{ex}[thm]{Example}

\newtheorem{rmk}[thm]{Remark}

\newcommand{\bproof}{\noindent{\textit{Proof. }}}
\newcommand{\eproof}{\hfill\qed}

\newcommand{\cm}{\mathcal{C}^{\mathcal{M}}}
\newcommand{\cmm}{\mathcal{C}^{\mathcal{M}}_m}
\newcommand{\cmzero}{\mathcal{C}^{\mathcal{M}}_0}

\newcommand{\cendmzero}{\mathcal{C}^{End_\Lambda(\mathcal{M})}_0}
\newcommand{\cendm}{\mathcal{C}^{End_\Lambda(\mathcal{M})}}
\newcommand{\cendmA}{\mathcal{C}^{End_\Lambda(\mathcal{M})}_A}

\newcommand{\aendmzero}{\mathcal{A}^{End_\Lambda(\mathcal{M})}_0}
\newcommand{\aendm}{\mathcal{A}^{End_\Lambda(\mathcal{M})}}
\newcommand{\aendmA}{\mathcal{A}^{End_\Lambda(\mathcal{M})}_A}
\newcommand{\am}{\mathcal{A}^{\mathcal{M}}}
\newcommand{\amm}{\mathcal{A}^{\mathcal{M}}_m}
\newcommand{\amzero}{\mathcal{A}^{\mathcal{M}}_0}

\DeclareFontFamily{U}{mathb}{\hyphenchar\font45}
\DeclareFontShape{U}{mathb}{m}{n}{
	<-6> mathb5 <6-7> mathb6 <7-8> mathb7
	<8-9> mathb8 <9-10> mathb9
	<10-12> mathb10 <12-> mathb12
}{}
\DeclareSymbolFont{mathb}{U}{mathb}{m}{n}
\DeclareMathSymbol{\lprec}{\mathrel}{mathb}{"CE}
\DeclareMathSymbol{\gprec}{\mathrel}{mathb}{"CF}

\begin{document}
	
	\title{Recursion relations and BPS-expansions in the HOMFLY-PT skein of the solid torus}
	\author{Lukas Nakamura}
	\maketitle
	
\begin{abstract}
	Inspired by the skein valued open Gromov-Witten theory of Ekholm and Shende and the Gopakumar-Vafa formula, we associate to each pair of non-negative integers $(g,l)$ a formal power series with values in the HOMFLY-PT skein of a disjoint union of $l$ solid tori. The formal power series can be thought of as open BPS-states of genus $g$ with $l$ boundary components and reduces to the contribution of a single BPS state of genus $g$ for $l=0$. Using skein theoretic methods we show that the formal power series satisfy gluing identities and multi-cover skein relations corresponding to an elliptic boundary node of the underlying curves. For $(g,l)=(0,1)$ we prove a crossing formula which is the multi-cover skein relation corresponding to a hyperbolic boundary node, also known as the pentagon identity.  
\end{abstract}

\section{Introduction}
This paper concerns formal power series in the HOMFLY-PT skein module of solid tori using relative recursion relations, motivated by skein valued curve counts \cite{es19} and the Gopakumar-Vafa formula \cite{gv98b}. We describe the background and motivation in Section \ref{ssec: motivation} and then state our results in Section \ref{ssec: main results}.  

\subsection{Background and motivation}\label{ssec: motivation}
Let $(M,\omega)$ be a symplectic Calabi-Yau threefold, i.e. $M$ is $6$-dimensional symplectic manifold and $\omega$ is a symplectic form on $M$ with vanishing associated first Chern class. Let $L \subset M$ be Lagrangian submanifold with vanishing Maslov class. Given some additional data, Ekholm and Shende \cite{es19} introduced, motivated by connections between skein theory, Chern-Simons theory, and open topological string theory \cite{wit89,wit95,ov00}, a skein valued open Gromov-Witten partition function $\GW(M,L)$ of $(M,L)$ that counts holomorphic curves in $M$ with boundary on $L$ by the values of their boundaries in the framed HOMFLY-PT skein of the $3$-manifold $L$, which is invariant under deformations. 

If $M$ is non-compact with a positive cylindrical end (convex at infinity) and if $L$ is cylindrical at infinity, then the boundaries of one-dimensional moduli spaces of holomorphic curves asymptotic to an index $1$ Reeb chord of the boundary $\partial L$ at infinity give relations in the relative skein of $L$. In basic cases, this leads to recursion relations that $\GW(M,L)$ satisfies, which were used in \cite{es20} to compute the Gromov-Witten partition function of a toric brane in $\C^3$ and of the union of the zero section in $T^*S^3$ with a push-off of the unknot conormal. In the case of a toric brane, the only embedded holomorphic curve is a disc, all other curves appear as multiple covers of this disc with constant curves attached. In the case of the unknot conormal, similarly, the only embedded holomorphic curve is an annulus between the push-off of the conormal and the zero-section and all other curves are multiple covers of this annulus with constant curves attached. 

In \cite{es21}, the latter example was generalized to a calculation of the local contribution to the Gromov-Witten partition function from a general basic annulus via the relations from boundaries of one-dimensional moduli spaces. Using this result, \cite{es21} further shows that the HOMFLY-PT generating function of any link in the $3$-sphere equals the Gromov-Witten partition function of a push-off of the link conormal, either in the resolved conifold or together with the zero-section in $T^*S^3$, thereby mathematically establishing the Ooguri-Vafa conjecture \cite{ov00}. 

Motivated by open Gromov-Witten theory and the Gopakumar-Vafa formula \cite{gv98a,gv98b}, we study in this article recursion relations for partition functions with values in the skein of a solid torus via skein-theoretic methods. We use recursion relations to study a certain family of skein valued partition functions which we call \emph{BPS partition functions}. For future reference, recall the Gopakumar-Vafa formula \cite{gv98b} which asserts that the closed Gromov-Witten partition function of a Calabi-Yau 3-fold can be expressed as 
\begin{equation}\label{eq:gopakumar-vafa formula}
	\GW(M, \emptyset) = \exp \left( \sum_{A \neq 0,g} \GW_{A,g}	x^{2g-2} q^A \right) = \prod_{A \neq 0, h} \exp \left(n_{A,h} \sum_{k=1}^\infty \frac{1}{k} \left(s^{k} - s^{-k} \right)^{2h-2} q^{kA} \right)
\end{equation}
for some integers $n_{A,h}$ which vanish for fixed $A$ and sufficiently large $h$. Here, the sum and the product are taken over all non-zero homology classes $A \in H_2(M;\Z)$ and all $g,h \geq 0$, and $\GW_{A,g}$ denotes the rational count of all holomorphic curves in $M$ of genus $g$ in the homology class $A$. The variables $x$ and $q$ are formal variables, and $s$ and $x$ are related via $s = \exp(i\frac{x}{2})$. The existence and uniqueness of the \emph{BPS numbers} $n_{A,h} \in \Q$ was proven by Bryan and Pandharipande in \cite{bp01}, their integrality was proven by Ionel and Parker \cite{ip18}, and their vanishing for large $h$ was proven by Doan, Ionel, and Walpuski \cite{diw21}. 

\subsection{Main results}\label{ssec: main results}
Motivated by the contribution of one BPS state on the right hand side of (\ref{eq:gopakumar-vafa formula}) (i.e., exactly one of the numbers $n_{A,h}$ is equal to one, and all the others are zero), we introduce for any numbers $g \in \Z$ and $l \in \Z_{\geq 0}$, which we refer to as \emph{the genus} and \emph{the number of boundary components} respectively, the partition functions
\begin{equation}\label{eq:def of bps part fct}
	\Psi^{\pm(g,l)}(t) \coloneqq \exp\left(\pm \sum\limits_{i=1}^\infty \frac{1}{i} \frac{t^i}{(s^i - s^{-i})^{2-2g-l}} (P_i)^{\otimes l}  \right)
\end{equation}
with values in the skein of a disjoint union of $l$ solid tori, where the $P_i$ are certain elements in the skein of a solid torus corresponding to power sums (see (\ref{eq:def of P_i}) below). 

\begin{rmk}
	Motivated by string theory, it is conjectured that from the HOMFLY-PT partition function of an $l$-component link (see Example \ref{ex:homflypt partition function} below) one can extract integer invariants, the so-called LMOV numbers $\widehat{N}_{(\lambda_1,...,\lambda_l),g,Q}$ \cite{ov00,lm01,lm02,lmv00}. 
	
	The partition functions defined in \eqref{eq:def of bps part fct} are precisely the partition functions that have exactly one non-trivial LMOV number in the fundamental representation which is equal to $\pm 1$, i.e.\, $\Psi^{\pm(g,l)}$ has $\widehat{N}_{(\square,...,\square),g,0} = \pm 1$ and all other LMOV numbers vanish.
\end{rmk}

In Section \ref{sec:bps partition functions}, we prove the following:

\begin{thm}\label{thm:properties of bps partition functions}
	$\quad$
	\begin{itemize}	
		\item[$\mathrm{(a)}$]\emph{(inverse)} $\Psi^{-(g,l)}$ is the inverse of $\Psi^{(g,l)}$: 
		\begin{equation*}
			\Psi^{-(g,l)} \Psi^{(g,l)} = 1.
		\end{equation*}
		\item[$\mathrm{(b)}$]\emph{(crossing a Lagrangian)} If one boundary component of $\Psi^{-(g,l)}$ for $l \geq 1$ is mapped to an unknot, then the resulting partition function is the difference of two BPS partition functions of the same genus and $l-1$ boundary components:
		\begin{equation*}
			\Psi^{\pm(g,l-1)}(a t) \Psi^{\mp(g,l-1)}(a^{-1} t) = \left( 1^{\otimes l-1} \otimes \mathcal{U} \right)(\Psi^{\pm(g,l)}).
		\end{equation*}
		\item[$\mathrm{(c)}$]\emph{(gluing)} The BPS partition functions satisfy
		\begin{equation*}
			\left(\sum\limits_\lambda 1^{\otimes l_1-1} \otimes  W_\lambda^* \otimes W_\lambda^* \otimes 1^{\otimes l_2-1} \right)\left(\Psi^{(-1)^{k_1}(g_1,l_1)} \otimes \Psi^{(-1)^{k_2}(g_2,l_2)} \right) = \Psi^{(-1)^{k_1+k_2}(g_1 + g_2,l_1 + l_2 - 2)}
		\end{equation*}
		and
		\begin{equation*}
			\left(\sum\limits_{i=0}^{\infty} 1^{\otimes l_1-2} \otimes  W_{(i)}^* \otimes W_{(i)}^* \right) \left(\Psi^{\pm(g,l)} \right) = \Psi^{\pm(g+1,l - 2)}.
		\end{equation*}
		\item[$\mathrm{(d)}$]\emph{(disk partition function)} In the case that $(g=0,l=1)$, the BPS partition functions are given by 
		\begin{equation*}
			\Psi^{\pm(0,1)} = 1 + \sum_{n\geq1} \sum_{\lambda \vdash n}\sum_\lambda (\pm1)^{n} W_\lambda \prod_{\square \in \lambda} \frac{s^{\pm c(\square)}}{s^{h(\square)} - s^{-h(\square)}}.
		\end{equation*}
		
		\item[$\mathrm{(e)}$] \emph{(annulus partition function)} In the case that $(g=0,l=2)$, the BPS partition functions are given by
		\begin{equation*}
			\Psi^{(0,2)} = 1 + \sum_{n\geq1} \sum_{\lambda \vdash n} W_\lambda \otimes W_\lambda.
		\end{equation*}
		and
		\begin{equation*}
			\Psi^{-(0,2)} = 1 + \sum_{n \geq 1} (-1)^n \sum_{\lambda \vdash n} W_\lambda \otimes W_{\lambda'}.
		\end{equation*}
		
		\item[$\mathrm{(f)}$] In the case that $(g=1,l=1)$, the BPS partition functions are given by 
		\begin{equation*}
			\Psi^{(1,1)} = 1 + z \sum_{i=n}^\infty C_n
		\end{equation*}
		and 
		\begin{equation*}
			\Psi^{-(1,1)} = 1 - z \sum_{n=1}^\infty \overline{C}_n.
		\end{equation*}
	\end{itemize}	
\end{thm}

\begin{rmk} (1) Some explanations are in order. All the partition functions take values in a suitable localization and completion of the framed HOMFLY-PT skein algebra of a disjoint union of solid tori, which as a module over the ring $\Q[a^{\pm 1},s^{\pm 1}, \frac{1}{s-s^{-1}}]$ is generated by isotopy classes of framed, oriented links modulo the local skein relations
\begin{equation}
	\begin{gathered}
		\begin{tikzpicture}
			\tikzstyle{smallnode}=[circle, inner sep=0mm, outer sep=0mm, minimum size=0.5mm, draw=black, fill=black];
			\begin{knot}[
				ignore endpoint intersections=false
				clip width=5,
				clip radius=8pt,
				only when rendering/.style={
				}
				]
				
				\strand [thick, ->] (0,0)
				to (1,1);
				\strand [thick, ->] (1,0)
				to (0,1);
				\node at (1.4,0.5) {$-$};
				\strand [thick, ->] (1.8,0)
				to (2.8,1);
				\strand [thick, ->] (2.8,0)
				to (1.8,1);
				\node at (3.2,0.5) {$=$};
				\node at (3.6,0.5) {$z$};
				\strand [thick, ->] (3.8,0)
				to [out=45,in=-45, looseness=1.8] (3.8,1);
				\strand [thick, ->] (4.8,0)
				to [out=135,in=225, looseness=1.8] (4.8,1);
				
				\flipcrossings{2}
			\end{knot}
		\end{tikzpicture}\\
		\begin{tikzpicture}
			\tikzstyle{smallnode}=[circle, inner sep=0mm, outer sep=0mm, minimum size=0.5mm, draw=black, fill=black];
			\begin{knot}[
				ignore endpoint intersections=false
				clip width=5,
				clip radius=8pt,
				only when rendering/.style={
				}
				]
				
				\strand [thick] (0,0)
				to [out=90, in=-90] (0,0.2)
				to [out=90, in=180] (0.2,0.7)
				to [out=0, in=90] (0.4,0.5);
				\strand [thick, ->] (0.4,0.5)
				to [out=-90, in=0] (0.2,0.3)
				to [out=180, in=-90] (0,0.8)
				to [out=90, in=-90] (0,1); 
				
				\node at (0.7,0.5) {$=$};
				\node at (1.1, 0.5) {$a$};
				
				\strand [thick, ->] (1.3,0)
				to (1.3,1); 
			\end{knot}
		\end{tikzpicture}\\
		\begin{tikzpicture}
			\draw (0,0) circle (0.5);
			\node at (0.9,0) {$=$};
			\node at (1.7,0) {$\frac{a-a^{-1}}{s-s^{-1}}$};
		\end{tikzpicture}
	\end{gathered}
\end{equation}
where $z \coloneqq s - s^{-1}$.

The symbols $P_i$, $C_n$, and $W_\lambda$ all denote different elements in the skein of a solid torus. While $P_i$ and $C_i$ admit nice representation in terms of knot diagrams (see \eqref{eq:def of P_i} and \eqref{eq:def of C_m} below), the explicit description of $W_\lambda$ is more complicated \cite{am98}. Two important facts are that (1) the $W_\lambda$'s from an eigenbasis of the positive part of the skein of a solid torus for a certain operator that will be crucial in our discussion of recursion relations and (2) there is an isomorphism from the positive part of the skein algebra of a solid torus to the algebra of symmetric functions which maps $W_\lambda$ to the Schur function $s_\lambda$ \cite{luk05}. Their dual basis elements are denoted by $W^*_\lambda$.

The map $\mathcal{U}$ is induced by embedding a solid torus into $S^3$ as a tubular neighbourhood of a zero-framed unknot. In (b), it is interpreted as a map from the skein of the solid torus to the ground ring. In Section \ref{sec:homflypt skein} we recall the necessary background from skein theory in more detail. 

In parts (a), (c), (d), (e), and (f), we omitted the formal variable $t$ from the notation.

	(2) Part (e) of Theorem \ref{thm:properties of bps partition functions} is equivalent to the well-known Cauchy identity for Schur functions, see for instance \cite{mac95}. Part (f) was previously proved by Morton and Manch\'{o}n \cite[Theorem 1, Theorem 2]{mm08}. In Section \ref{sec:bps partition functions} we give new proofs of these results using recursion relations.
\end{rmk}

Geometrically, (a) can be interpreted as a wall-crossing identity of the moduli spaces under deformation: a BPS state and its anti-state are annihilated or created. 

Also (b) has a natural interpretation in skein valued Gromov-Witten theory: given a transversally cut-out embedded holomorphic curve $C$ in $M$ with boundary on $L$, consider a one parameter deformation where an interior point $p$ of $C$ crosses the Lagrangian at some instant $t_0$. As was observed in \cite{es19}, the difference of the curve at time $t = t_0 + \eps$ and the curve at time $t_0 - \eps$ equals the curve $C'$ obtained from $C$ by removing an open disk around $p$ and with the resulting boundary circle mapped to a small unknot near the intersection point in the Lagrangian. Here (b) asserts that this picture is valid not only for the underlying embedded curve but for the whole partition function it carries. 

In (d) and (e), $\Psi^{\pm(0,1)}$ and $\Psi^{(0,2)}$ are the Gromov-Witten partition functions of an isolated basic disk \cite{es20} and an isolated basic annulus \cite{es21}, respectively. Then, when $(g,l)=(0,2)$, (b) together with (c) and (d) imply that the HOMFLY-PT generating function of the unknot equals the product of two disk partition functions. 

The first identity in (c) describes the gluing of two BPS partition functions across their boundary tori, the second one the self-gluing of one partition function across two of its boundary tori. From the point of view of holomorphic curves such gluing identities appear via so called Symplectic Field Theory (SFT) stretching \cite{egh10}. In the general case, it works as follows. Consider two holomorphic curves $u_1$ and $u_2$ in $M$ with boundary on a Lagrangian $L$ of topology $S^{1}\times\R^{2}$. Consider a metric on $L$ with a unique simple geodesic loop $S^{1}\times\{0\}$. The corresponding contact form on the unit conormal bundle of $ST^{\ast}(S^{1}\times\R^{2})$ then has exactly two simple Reeb orbits $\gamma^{\pm1}$ corresponding to the geodesic parametrized with opposite orientations. Applying SFT-stretching along the unit conormal with respect to this contact form, $u_{j}$ decomposes into an outside curve $u_{j}'$ with negative punctures asymptotic to Reeb orbits $\gamma^{(\pm d)}$ that are multiples of the underlying simple Reeb orbits and an inside curve $u''_{j}$ with positive asymptotics at $\gamma^{(\pm d)}$, matching the negative asymptotics of $u'_j$, and boundaries on the zero-section $L$. Note that these inside curves are universal and determined by their positive punctures only.
Applying the gluing operation to the curves $u_{1}$ and $u_{2}$ then means joining $u'_1$ and $u'_2$ by gluing them with all possible curves in $T^{\ast} L$ that join their asymptotics.

In (c), we consider the case when the boundaries of $u_+ \coloneqq u_{1}$ and $u_- \coloneqq u_{2}$ parametrize the simple geodesic with opposite orientations. Then considering also multiple covers of the underlying simple curve, the gluing formula is compatible with the following picture: after SFT-stretching, $u_{\pm}$ and its multiple covers break into a simple curve $u'_\pm$ on the outside with one negative puncture asymptotic to $\gamma^{(\pm1)}$ and its multiple covers and a one-punctured holomorphic disk $u^d_\pm$ on the inside, which matches the asymptotic of $u'_\pm$ and has boundary on $L$, and its multiple covers whose total skein-valued Gromov-Witten count is given by
\begin{equation}
	u^d_\pm = \exp\left( \sum_{i=1}^{\infty} \frac{1}{i} \gamma^{(\pm i)} \otimes P_i \right).
\end{equation}
Then the first identity in (c) expresses that the Gromov-Witten count of 
\begin{equation}
	u \coloneqq u_+ * u^s * u_-
\end{equation}
is equal to $\Psi^{(-1)^{k_1+k_2}(g_1 + g_2,l_1 + l_2 - 2)}$ if the counts of $u_+$ and $u_-$ are equal to $\Psi^{(-1)^{k_1}(g_1,l_1)}$ and $\Psi^{(-1)^{k_2}(g_2,l_2)}$, where $*$ denotes the SFT-gluing of holomorphic curves and their multiple covers and $u^s$ denotes a two-punctured holomorphic sphere on the inside with positive asymptotics at $\gamma$ and $\gamma^{(-1)}$ whose total count is given by
\begin{equation}
	u^s = \exp\left( \sum_{i=1}^{\infty} \frac{1}{i} \gamma^{(i)} \otimes \gamma^{(-i)} \right).
\end{equation}

As $\Psi^{(g,l)}$, $u^d_\pm$, and $u^s$ have natural expressions as the exponential of a partition function, the gluing formula for the BPS partition functions can also be interpreted as the exponentiated version of SFT-gluing of \emph{connected} curves. \\

Finally, we consider the hyperbolic counterpart of (b) in the simplest case. In \cite{es19}, it was observed that in a generic one-parameter family, when the boundaries of holomorphic curves intersect the difference of their Gromov-Witten contributions equals the contribution of the curve obtained by gluing at the intersection point. For basic disks, this is true not only for the underlying embedding, but for the whole partition function (see Theorem \ref{thm:Psi^(0,1) crossing formula} for a precise statement):

\begin{thm}\label{thm:multicover disk crossing} \textnormal{(multicover disk crossing)} Two disk partition  functions before crossing are equal to two disk partition functions after crossing together with the partition function of the glued disk.
\end{thm}

If we restrict to the $U(1)$-skein and embed the boundary connected sum of the two solid tori corresponding to the boundaries of the two disk partition functions into a single solid torus, then Theorem \ref{thm:multicover disk crossing} reduces to the multi-cover skein relation of Ekholm, Kucharski, and Longhi \cite{ekl20}.\\

Parts (a) and (b) of Theorem \ref{thm:properties of bps partition functions} are easy consequence of the definition of $\Psi^{\pm(g,l)}$. To prove Theorem \ref{thm:multicover disk crossing} and the other parts of Theorem \ref{thm:properties of bps partition functions}, we use relative recursion relations for partition functions with values in the skein of a disjoint union of solid tori which are introduced in Section \ref{sec:skein valued recursion relations}.

\begin{rmk}
	During the preparation of this article, Hu, Schrader, and Zaslow \cite{hsz23} and Scharitzer and Shende \cite{ss23} independently proved weaker versions of the disk crossing formula, called \emph{pentagon identity} in their work. We were also informed that Mingyuan Hu \cite{hu24} independently obtained a proof of Theorem \ref{thm:multicover disk crossing}.
\end{rmk}

\subsection{Overview} 
In Section \ref{sec:homflypt skein}, we recall definitions from skein theory and give an overview over the different versions of the HOMFLY-PT skeins we consider in this article. Section \ref{sec:skein valued recursion relations} develops the theory of (relative) recursion relations for partition functions with values in the skein of a solid torus. The main result of that section is an existence result for solutions of relative recursion relations (Theorem \ref{thm:relative recursion relation with solution}). In Section \ref{sec:bps partition functions}, we introduce the BPS partition functions $\Psi^{\pm(g,l)}$ and prove Theorem \ref{thm:properties of bps partition functions}.
In Section \ref{sec:multicover disk crossings}, we prove the skein valued multicover disk crossing formula.\\

\subsection*{Acknowledgements} I am very grateful to my advisor Tobias Ekholm for many helpful discussions about the contents of this paper. This work is partially supported by the Knut and Alice Wallenberg Foundation, grant KAW2020.0307.

\section{The HOMFLY-PT skein}\label{sec:homflypt skein}

In this section, we recall necessary definitions from skein theory and describe the geometric setup for the subsequent sections.

Let $\Lambda$ denote the the polynomial ring $\Q[a^\pm,s^\pm]$ localized\footnote{To define the framed HOMFLY-PT skein, it is not necessary to localize $\Q[a^\pm,s^\pm]$, but it turns out to be convenient for the study of the skein of the solid torus.} to allow for denominators of the forms $z \coloneqq s - s^{-1}$ and $\sum_{(i,j) \in \lambda} s^{2(j-i)}$ for non-empty partitions $\lambda \vdash n\in\N$, where $(i,j)$ denotes the coordinates in the Young diagram of $\lambda$. Given a $3$-manifold $M$ with a choice of \emph{in-going} and \emph{out-going} boundary points, we define the framed HOMFLY-PT skein $\mathcal{S}(M)$ of $M$ to be free $\Lambda$-module generated by finite disjoint collections of oriented framed arcs connecting the in-going to the out-going boundary points and oriented framed links in $M$ modulo framed smooth isotopy and the HOMFLY-PT skein relations 
\begin{equation}
	\begin{gathered}
		\begin{tikzpicture}
			\tikzstyle{smallnode}=[circle, inner sep=0mm, outer sep=0mm, minimum size=0.5mm, draw=black, fill=black];
			\begin{knot}[
				ignore endpoint intersections=false
				clip width=5,
				clip radius=8pt,
				only when rendering/.style={
				}
				]
				
				\strand [thick, ->] (0,0)
				to (1,1);
				\strand [thick, ->] (1,0)
				to (0,1);
				\node at (1.4,0.5) {$-$};
				\strand [thick, ->] (1.8,0)
				to (2.8,1);
				\strand [thick, ->] (2.8,0)
				to (1.8,1);
				\node at (3.2,0.5) {$=$};
				\node at (3.6,0.5) {$z$};
				\strand [thick, ->] (3.8,0)
				to [out=45,in=-45, looseness=1.8] (3.8,1);
				\strand [thick, ->] (4.8,0)
				to [out=135,in=225, looseness=1.8] (4.8,1);
				
				\flipcrossings{2}
			\end{knot}
		\end{tikzpicture}\\
		\begin{tikzpicture}
			\tikzstyle{smallnode}=[circle, inner sep=0mm, outer sep=0mm, minimum size=0.5mm, draw=black, fill=black];
			\begin{knot}[
				ignore endpoint intersections=false
				clip width=5,
				clip radius=8pt,
				only when rendering/.style={
				}
				]
				
				\strand [thick] (0,0)
				to [out=90, in=-90] (0,0.2)
				to [out=90, in=180] (0.2,0.7)
				to [out=0, in=90] (0.4,0.5);
				\strand [thick, ->] (0.4,0.5)
				to [out=-90, in=0] (0.2,0.3)
				to [out=180, in=-90] (0,0.8)
				to [out=90, in=-90] (0,1); 
				
				\node at (0.7,0.5) {$=$};
				\node at (1.1, 0.5) {$a$};
				
				\strand [thick, ->] (1.3,0)
				to (1.3,1); 
			\end{knot}
		\end{tikzpicture}\\
		\begin{tikzpicture}
			\draw (0,0) circle (0.5);
			\node at (0.9,0) {$=$};
			\node at (1.7,0) {$\frac{a-a^{-1}}{s-s^{-1}}$};
		\end{tikzpicture}
	\end{gathered}
\end{equation}
where we use the blackboard framing to (locally) identify a framed link with its projection onto the plane. Here, we also require that the boundary points come equipped with a choice of vector tangent to $\partial M$ and that the framing of the arcs agrees with the direction of that vector at the boundary point.

\subsection{Wiring} Given a smooth embedding of $M$ into another $3$-manifold $N$ together with a choice of a framed tangle in $\overline{N \setminus M}$ connecting the in-going boundary points of $\overline{N \setminus M}$ to the out-going boundary points of $\overline{N \setminus M}$ where we view the out/in-going boundary points of $M$ as additional in/out-going boundary points of $\overline{N \setminus M}$, we obtain an induced map $\mathcal{S}(M) \to \mathcal{S}(N)$ given by mapping a linear combination of framed tangles in $M$ to the corresponding linear combination of the union of the framed tangles in $M$ with the framed tangle in $N$. More generally, also a finite formal linear combinations $W$ of such tangles with coefficients in $\Lambda$ induces a map $W:\mathcal{S}(M) \to \mathcal{S}(N)$, which we will denote by the same symbol. $W$ is called a \emph{wiring}.

\subsection{The skein of a surface}
Given a $2$-dimensional surface $\Sigma$ with a choice of in-going and out-going boundary points, we write $\mathcal{S}(\Sigma) \coloneqq \mathcal{S}(\Sigma \times [0,1])$. Here, an in/out-going boundary point $p \in \partial \Sigma$ becomes the in/out-going boundary point $p \times \frac{1}{2}$ of $\Sigma \times [0,1]$, and the vector at $p \times \frac{1}{2}$ is given by $\partial_t$, where $t$ denotes the coordinate on $[0,1]$. We call $\mathcal{S}(\Sigma)$ \emph{the skein of $\Sigma$}.

If there are no boundary points, $\mathcal{S}(\Sigma)$ has the structure of a $\Lambda$-algebra with multiplication
\begin{equation}
	\mathcal{S}(\Sigma) \otimes \mathcal{S}(\Sigma) \cong \mathcal{S}\left( \Sigma \coprod \Sigma \right) \to \mathcal{S}(\Sigma)
\end{equation}
which is induced by the two embeddings of $\Sigma \times [0,1]$ into $\Sigma \times [0,1]$ coming from two orientation preserving embeddings of $[0,1]$ into $[0,1]$ as $[0,\frac{1}{2}]$ or $[\frac{1}{2},1]$, respectively.

There is a $\Q$-linear automorphism of $\mathcal{S}(\Sigma)$, called the \emph{mirror map}, induced by reversing all crossings of a link diagram and sending $s$ to $s^{-1}$ and $a$ to $a^{-1}$.

\subsection{The skein of $S^3$}

The skein $\mathcal{S}(S^3)$ is isomorphic to $\Lambda$ (see \cite{fyhlmo85,pt87}), i.e. any linear combinations of framed links in $S^3$ can be reduced in a unique way to a multiple of the empty knot via framed isotopy and the skein relations. Note that the skein of any $3$-manifold is a module over $\mathcal{S}(S^3)$, where the module structure is induced by the connected sum with $S^3$. Under the above isomorphism $\mathcal{S}(S^3) \approx \Lambda$, this module structure agrees with the structure as a $\Lambda$-module.

\subsection{The skein of a square with $2n$ boundary points}\label{sec:R^n_n} Denote by $R^n_n$ the skein of the square $[0,1] \times [0,1]$ with $n$ in-going boundary points on the bottom of the square and $n$ out-going boundary points on the top of the square for any $n \in \N_0$. Stacking two such squares on top of each other induces a map $R^n_n \otimes R^n_n \to R^n_n$ which induces the structure of a $\Lambda$-algebra on $R^n_n$. This algebra is isomorphic to the Hecke algebra $H_n$ of type $A$ of the symmetric group as was shown in \cite[Section 7]{mt90}.\\

Aiston and Morton \cite{am98} defined a set $e_\lambda \in R^n_n$ of mutually orthogonal quasi-idempotents indexed by partitions $\lambda$ of $n$, which after normalizing give rise to mutually orthogonal idempotents $E_\lambda$. These elements behave well under multiplication by the following two central elements in $R^n_n$: a \emph{full twist} 
\begin{equation}\label{eq:definition of F}
	\begin{tikzpicture}
		\tikzstyle{smallnode}=[circle, inner sep=0mm, outer sep=0mm, minimum size=0.5mm, draw=black, fill=black];
		\begin{knot}[
			ignore endpoint intersections=false
			clip width=5,
			clip radius=8pt,
			only when rendering/.style={
			}
			]
			
			\node at (-1.2,1) {$F$};
			\node at (-0.8,1) {$\coloneqq$};
			\strand [thick] (0,0)
			to [out=90, in=-90] (0,0.9)
			to [out=90, in=180] (0.1,1.1)
			to [out=0, in=90] (0.2,1)
			to [out=-90, in=0] (0.1,0.9);
			\strand [thick, ->] (0,1.55)
			to [out=90, in=-90] (0,2); 
			
			\strand [thick] (-0.4,0)
			to [out=90, in=-90] (-0.4,0.5)
			to [out=90, in=180, looseness=1.2] (0.1,1.5)
			to [out=0, in=90] (0.6,1)
			to [out=-90, in=0] (0.1,0.5)
			to [out=180, in=-10] (0.05,0.52);
			\strand [thick, ->] (-0.38,1.15)
			to [out=95, in=-90] (-0.4,1.4)
			to [out=90, in=-90] (-0.4,2); 
			
			\node at (-0.2,0.2) {$...$};
			\node at (-0.2,1.8) {$...$};
			\node at (0.4,1) {$...$};
			
			\flipcrossings{}
		\end{knot}
	\end{tikzpicture}
\end{equation}
and a meridian encircling the identity diagram in $R^n_n$
\begin{equation}\label{eq:definition of P}
	\begin{tikzpicture}
		\tikzstyle{smallnode}=[circle, inner sep=0mm, outer sep=0mm, minimum size=0.5mm, draw=black, fill=black];
		\begin{knot}[
			ignore endpoint intersections=false
			clip width=5,
			clip radius=8pt,
			only when rendering/.style={
			}
			]
			
			\node at (-1.4,1) {$P$};
			\node at (-1,1) {$\coloneqq$};
			\strand [thick, ->] (0,0)
			to [out=90, in=-90] (0,2); 
			
			\strand [thick, ->] (-0.4,0)
			to [out=90, in=-90] (-0.4,2); 
			
			\strand [thick, ->] (-0.44,1.12)
			to [out=190,in=90] (-0.6,1)
			to [out=-90,in=180] (-0.2,0.8)
			to [out=0,in=-90] (0.2,1)
			to [out=90,in=-10] (0.04,1.12);
			
			\node at (-0.2,0.2) {$...$};
			\node at (-0.2,1.8) {$...$};
			
			\flipcrossings{1,2}
		\end{knot}
	\end{tikzpicture}.
\end{equation}
Namely, Aiston and Morton showed that $F e_\lambda = a^{|\lambda|} s^{\sum_{(i,j) \in \lambda}2(j-i)} e_\lambda$ and that 
\begin{equation}\label{eq:p-unknot on e_lambda}
(P - \bigcirc) e_\lambda = a z \sum_{(i,j) \in \lambda} s^{2(j-i)} e_\lambda,
\end{equation}
where $\bigcirc$ denotes multiplication by $\frac{a-a^{-1}}{s-s^{-1}}$.

Blanchet \cite{bla00} showed that $R^n_n$ has a basis $\alpha_t \beta_\tau$ indexed by pairs of standard tableaus $(t,\tau)$ for partitions $\lambda= \lambda(t) = \lambda(\tau) \vdash n$, where each $\alpha_t \beta_\tau$ lies in the two-sided ideal generated by $e_{\lambda}$. In particular, it follows that
\begin{equation}\label{eq:action of p-unknot on R^n_n}
	(P - \bigcirc) \left( \alpha_t \beta_\tau \right) = a z \sum_{(i,j) \in \lambda} s^{2(j-i)} \alpha_t \beta_\tau.
\end{equation}
Thus, $P-\bigcirc: R^n_n \to R^n_n$ is invertible if $n \geq 1$ since $a z \sum_{(i,j) \in \lambda} s^{2(j-i)} \in \Lambda$ is a unit. In fact, even more is true: let $M$ be a $3$-manifold with boundary, and let $B \times [0,1] \subseteq M$ be an embedded cylinder, where $B$ denotes a 2-disk, so that $\partial B \times [0,1] \subseteq \partial M$ does not meet any of the chosen boundary points of $M$. Choose a homeomorphism from the cube $[0,1]^3$ to $B \times [0,1]$ which restricts to a diffeomorphism on the complement of the vertices and edges of $[0,1]^3$. This induces a choice of $n$ additional out-going boundary points in $M \setminus B \times (0,1)$ on $B \times \{0\}$ and a choice of $n$ additional in-going boundary points on $B \times \{1\}$ from the above choice of boundary points on $[0,1]^3$. We obtain a decomposition of $M$ into $[0,1]^3$ and $M \setminus B \times (0,1)$ compatible with the chosen boundary points which induces a map $W: R^n_n \otimes \mathcal{S}(M \setminus B \times (0,1)) \to \mathcal{S}(M)$. Denote by $P: \mathcal{S}(M) \to \mathcal{S}(M)$ the linear map induced by encircling $B \times \{\frac{1}{2}\}$ by its boundary. It follows from the construction that $W \circ (P \otimes 1) = P \circ W$. Using this and (\ref{eq:action of p-unknot on R^n_n}), it follows that the restriction of $P - \bigcirc$ to the image of $W$ is invertible if $n \geq 1$ and we assume that $\mathcal{S}(M)$ is a free $\Lambda$-module. Examples of 3-manifolds $M$ whose skein is a free $\Lambda$-module are handlebodies \cite{prz92} and the skein of a solid torus with two boundary points (see Section \ref{sec:skein of torus with two boundary points} below).

\subsection{The skein $\mathcal{C}$ of the solid torus} Let $\mathcal{C} \coloneqq \mathcal{S}(S^1 \times [0,1])$ denote the skein of a solid torus. Note that the map $\mathcal{C} \otimes \mathcal{C} \to \mathcal{C}$ given by the embedding the disjoint union of two annuli into a bigger annulus induces the structure of a commutative $\Lambda$-algebra on $\mathcal{C}$.

Turaev \cite[Theorem 2]{tur90} showed that $\mathcal{C}$ is isomorphic (as an algebra) to the polynomial algebra $\Lambda[C_m]$ in countably infinitely many variables $C_m, m\in \Z_{\neq0},$ which are defined as 
\begin{equation}\label{eq:def of C_m}
	\begin{tikzpicture}
		\tikzstyle{smallnode}=[circle, inner sep=0mm, outer sep=0mm, minimum size=0.5mm, draw=black, fill=black];
		\begin{knot}[
			ignore endpoint intersections=false
			clip width=5,
			clip radius=8pt,
			only when rendering/.style={
			}
			]
			\node at (-1,0.5) {$C_m$};
			\node at (-0.4,0.5) {$=$};
			
			\strand [thick] (0,0)
			to (1.4,1);
			\strand [thick] (0.2,0)
			to (0,1);
			\strand [thick] (0.4,0)
			to (0.2,1);
			\strand [thick] (0.601,0)
			to (0.4,1);
			\strand [thick] (1.4,0)
			to (1.2,1);
			
			\node at (1,0.3) {$...$};
			
		\end{knot}
	\end{tikzpicture},
\end{equation}
where we identify the top and the bottom to get an element in $\mathcal{C}$ and $m$ denotes the number of strands in the figure which are oriented from bottom to top if $m>0$ and from top to bottom if $m<0$.\\

We denote by $\widehat{(\cdot)}: R^n_n \to \mathcal{C}$ the map induced by mapping a tangle $T$ to its closure in the solid torus which is obtained by identifying the top with the bottom of the cube. Denote its image by $\mathcal{C}_n$. Then 
\begin{equation}
\mathcal{C}^+ \coloneqq \bigoplus_{n \in \N_0} \mathcal{C}_n
\end{equation}
is an $\N_0$-graded $\Lambda$-algebra with unit $1 \in \mathcal{C}_0$, where the degree-$n$ part is given by $\mathcal{C}_n$. Note that $\mathcal{C}^+ = \Lambda[C_1, C_2,...]$ as a $\Lambda$-algebra and that $C_n \in \mathcal{C}_n$ for $n > 0$.

Let $P: \mathcal{C} \to \mathcal{C}$ be given by encircling the solid torus by a simple closed curve analogously to how $P$ is defined on $R^n_n$ (see (\ref{eq:definition of P}) above). $P$ maps $\mathcal{C}_n$ into $\mathcal{C}_n$ since $\widehat{(\cdot)}\circ P = P \circ \widehat{(\cdot)}$, and $P|_{\mathcal{C}^+}$ has an eigenbasis with $1$-dimensional eigenspaces spanned by the closures $W_\lambda \coloneqq \widehat{E_\lambda}$ of the idempotents in $R^n_n$. As shown above, $(P - \bigcirc)|_{\mathcal{C}^+_{>0}}$ is invertible. Note that $P - \bigcirc|_{\mathcal{C}_0} = 0$.

Similarly, the full twist $F: \mathcal{C} \to \mathcal{C}$ is defined as in (\ref{eq:definition of F}), and it acts on $W_\lambda$ as
\begin{equation}\label{eq:action of F on W_lambda}
	F \left(W_\lambda \right)= a^{|\lambda|} s^{\sum_{(i,j) \in \lambda}2(j-i)} W_\lambda.
\end{equation}

There exists an isomorphism of $\Lambda$-algebras from the ring of symmetric functions with coefficients in $\Lambda$ to $\mathcal{C}^+$ which maps the Schur function $s_\lambda$ to $W_\lambda$ (see \cite{luk05}). The elements $W_\lambda$ are invariant under the mirror map. In the case that $\lambda = (n)$ for $n \in \N$, this is shown in \cite{mor02}, and the general case follows immediately since any $W_\lambda$ can be written as a polynomial in the elements $W_{(n)}$ with coefficients in $\Q$.\\

A framed knot $K \subseteq S^3$ induces an embedding of the solid torus into $S^3$ which is well-defined up to isotopy. Denote the induced map by $K: \mathcal{C} \to \mathcal{S}(S^3)$ as well, i.e. for $A \in \mathcal{C}$, $K(A)$ is the satellite link with companion $K$ and pattern $A$. \emph{The $\lambda$-colored HOMFLY-PT polynomial of $K$} is defined as $K(W_\lambda) \in \mathcal{S}(S^3) = \Lambda$. More generally, consider a framed $l$-component link $L \subseteq S^3$ with components $L_1,...,L_l$. Then $L$ induces a map $L: \mathcal{C}^{\otimes l} \to \mathcal{S}(S^3)$, and we call $L(W_{\lambda_1} \otimes ... \otimes W_{\lambda_l}) \in \Lambda$ the \emph{$(\lambda_1,...,\lambda_l)$-colored HOMFLY-PT polynomial of $L$} for any partitions $\lambda_1,...,\lambda_l$.

\subsection{The skein of the annulus with two boundary points}\label{sec:skein of torus with two boundary points} Let $\mathcal{A} \coloneqq \mathcal{S}(S^1 \times [0,1])$ be the skein of the annulus with one in-going and one out-going boundary point. Up to $\Lambda$-linear isomorphism, $\mathcal{A}$ does not depend on which two boundary points we choose as the boundary of a solid torus is connected. Following \cite{mor02}, we draw the in-going boundary point on the left boundary of the annulus and the out-going point on the right. Then $\mathcal{A}$ naturally has the structure of a $\Lambda$-algebra with multiplication induced by the wiring
\begin{equation}
	\begin{tikzpicture}
		\tikzstyle{smallnode}=[circle, inner sep=0mm, outer sep=0mm, minimum size=0.7mm, draw=black, fill=black];
		\begin{knot}[
			draft mode=strands,
			ignore endpoint intersections=true
			clip width=5,
			clip radius=8pt,
			only when rendering/.style={
			},
			draft mode=strands
			]
			\node[smallnode] (left)   at (-0.5,0.75) {};
			\strand [thick] (-0.5,0.75)
			to (0,0.75);
			
			\strand [thick] (0,0)
			to (0,1.5);
			\strand [thick] (0.5,0)
			to (0.5,1.5);
			
			\strand [thick, ->] (0.5,0.75)
			to (1,0.75);
			\strand [thick] (1,0.75)
			to (1.5,0.75);

			\strand [thick] (1.5,0)
			to (1.5,1.5);
			\strand [thick] (2,0)
			to (2,1.5);
			
			\strand [thick] (2,0.75)
			to (2.5,0.75);	
			\node[smallnode] (left)   at (2.5,0.75) {};
		\end{knot}
	\end{tikzpicture}.
\end{equation}
Here, the areas between the two pairs of vertical lines represent two annuli (oriented from bottom to top), and the horizontal arcs describe the arcs connecting the boundary points. As usual, the top and the bottom are identified, and we do not draw the left and right boundary of the bigger annulus.\\

As shown in \cite{mor02}, the algebra $\mathcal{A}$ is commutative. In fact, $\mathcal{A}$ is also a $\mathcal{C}$-bimodule, and the multiplications of $\mathcal{C}$ and $\mathcal{A}$ are compatible under this structure. Indeed, define the $\Lambda$-algebra morphisms $l: \mathcal{C} \otimes \mathcal{A} \to \mathcal{A}$ and $r:\mathcal{A} \otimes \mathcal{C} \to \mathcal{A}$ by placing an annulus underneath or, respectively, above the annulus with the boundary points:
\begin{equation}
	\begin{tikzpicture}
		\tikzstyle{smallnode}=[circle, inner sep=0mm, outer sep=0mm, minimum size=0.7mm, draw=black, fill=black];
		\begin{knot}[
			ignore endpoint intersections=true,
			clip width=5,
			clip radius=8pt,
			only when rendering/.style={
			},
			]
			\node at (-1.1,0.75) {$l=$};
			
			\node[smallnode] (left)   at (-0.5,0.75) {};
			\strand [thick] (-0.5,0.75)
			to (0,0.75);
			
			\strand [thick] (0,0)
			to (0,1.5);
			\strand [thick] (0.5,0)
			to (0.5,1.5);
			
			\strand [thick] (0.5,0.75)
			to (2,0.75);
			\node[smallnode] (left)   at (2,0.75) {};
			
			\strand [thick] (1,0)
			to (1,1.5);
			\strand [thick] (1.5,0)
			to (1.5,1.5);

			\node at (2.2,0.75) {$,$};
			\node at (3,0.75) {$r=$};
			
			\node[smallnode] (left)   at (3.6,0.75) {};
			\strand [thick] (3.6,0.75)
			to (4.1,0.75);
			
			\strand [thick] (4.1,0)
			to (4.1,1.5);
			\strand [thick] (4.6,0)
			to (4.6,1.5);
			
			\strand [thick] (4.6,0.75)
			to (5.0,0.75);
			\strand [thick] (5.7,0.75)
			to (6.1,0.75);
			\node[smallnode] (left)   at (6.1,0.75) {};
			
			\strand [thick] (5.1,0)
			to (5.1,1.5);
			\strand [thick] (5.6,0)
			to (5.6,1.5);
		\end{knot}
	\end{tikzpicture}.
\end{equation} 
We denote the class of the simple arc from the in-going to the out-going point by $e$, and we denote the arc going $m \in \Z$ times around the annulus from the in-going to the out-going boundary point by $c^m$. Clearly, $c^0 = e$ is the unit in $\mathcal{A}$, and $c^m c^n = c^{m+n}$.
\begin{equation}
	\begin{tikzpicture}
		\tikzstyle{smallnode}=[circle, inner sep=0mm, outer sep=0mm, minimum size=0.7mm, draw=black, fill=black];
		\begin{knot}[
			ignore endpoint intersections=true,
			clip width=5,
			clip radius=8pt,
			only when rendering/.style={
			},
			]
			\node at (-0.6,0.75) {$c^m=$};
			\node[smallnode] at (0,0.75) {};
			\strand [thick] (0,0.75)
			to [out=0,in=-90] (0.3,1.5);
			\strand [thick] (0.3,0)
			to [out=90,in=-90] (0.4,1.5);
			\strand [thick] (0.4,0)
			to [out=90,in=-90] (0.5,1.5);
			\node at (0.7,0.75) {$...$};
			\strand [thick] (0.9,0)
			to [out=90,in=-90] (1,1.5);
			\strand [thick] (1,0)
			to [out=90,in=180] (1.3,0.75);
			\node[smallnode] at (1.3,0.75) {};
			\node at (1.5,0.75) {$,$};
			
			\node at (2.38,0.75) {$e=$};
			\node[smallnode] at (2.9,0.75) {};
			\strand [thick] (2.9,0.75)
			to [out=0,in=180] (4.2,0.75);
			\node[smallnode] at (4.2,0.75) {};
		\end{knot}
	\end{tikzpicture}.
\end{equation}
We also simply write $c \coloneqq c^1$. Using the unit $e$, $l$ and $r$ induce morphisms of $\Lambda$-algebras $l(\cdot,e), r(e,\cdot): \mathcal{A} \to \mathcal{C}$. We denote their difference by $[\cdot,e] \coloneqq l(\cdot,e) - r(e,\cdot)$.\\

For $n \geq 0$, let $W_n$ be the wiring
\begin{equation}\label{eq:wiring W_n}
\begin{tikzpicture}
\tikzstyle{smallnode}=[circle, inner sep=0mm, outer sep=0mm, minimum size=0.7mm, draw=black, fill=black];
\begin{knot}[
ignore endpoint intersections=false
clip width=5,
clip radius=8pt,
flip crossing=1,
only when rendering/.style={
}
]
\node at (-2.5,1.6) {$W_n =$};
\draw[draw=black] (-1.5,2) rectangle ++(1.4,0.4);
\node[smallnode] (left)   at (-1.8,1.6) {};
\strand [thick] (-1.8,1.6)
to [out=right, in=down, looseness=1.5] (-1.3,2);
\node at (-1.1,0.3) {$n$};
\strand [thick] (-1.1,0.6)
to [out=up, in=down] (-1.1,2);
\strand [thick] (-0.9,0.6)
to [out=up, in=down] (-0.9,2);
\strand [thick] (-0.7,0.6)
to [out=up, in=down] (-0.7,2);
\strand [thick] (-0.5,0.6)
to [out=up, in=down] (-0.5,2);
\strand [thick] (-0.3,0.6) node[below]{$1$}
to [out=up, in=down] (-0.3,2);
\strand [thick,-{>[scale=1.7,length=2,width=2]}] (-1.3,2.4)
to [out=up, in=down, looseness=1.5] (-1.1,2.8);
\strand [thick,-{>[scale=1.7,length=2,width=2]}] (-1.1,2.4)
to [out=up, in=down, looseness=1.5] (-0.9,2.8);
\strand [thick,-{>[scale=1.7,length=2,width=2]}] (-0.9,2.4)
to [out=up, in=down, looseness=1.5] (-0.7,2.8);
\strand [thick,-{>[scale=1.7,length=2,width=2]}] (-0.7,2.4)
to [out=up, in=down, looseness=1.5] (-0.5,2.8);
\strand [thick,-{>[scale=1.7,length=2,width=2]}] (-0.5,2.4)
to [out=up, in=down, looseness=1.5] (-0.3,2.8);
\node[smallnode] (right)   at (0.9,1.6) {};
\strand [thick] (-0.3,2.4)
to [out=up, in=up, looseness=1.5] (0.1,2.4)
to [out=down, in = up] (0.1,1.8)
to [out=down, in = left, looseness = 0.8] (0.9,1.6);

\end{knot}
\end{tikzpicture}
\end{equation}
We denote by $\mathcal{A}_n \coloneqq W_n(R^{n+1}_{n+1})$ the image of the induced map $W_n: R^{n+1}_{n+1} \to \mathcal{A}$, and we define $\mathcal{A}^+ \coloneqq \bigoplus_{n\in\N_0} \mathcal{A}_n$. As for $\mathcal{C}^+$, this decomposition induces an $\N_0$-grading on $\mathcal{A}^+$. Lukac \cite[Lemma 3.3.3]{luk01} showed that as a graded algebra over the graded ring $\mathcal{C}^+$, $\mathcal{A}^+$ is the polynomial algebra $\mathcal{C}^+[c]$ with respect to both algebra structures $l$ and $r$. In particular, $\mathcal{A}^+$ is a free $\Lambda$-module.\\

We denote by $\text{cp}:\mathcal{A} \to \mathcal{C}$ the map induced by connecting the out-going boundary point to the in-going boundary point by an arc going underneath the solid torus: 
\begin{equation}\label{eq:definition of cp}
	\begin{tikzpicture}
		\tikzstyle{smallnode}=[circle, inner sep=0mm, outer sep=0mm, minimum size=0.7mm, draw=black, fill=black];
		\begin{knot}[
			ignore endpoint intersections=true,
			clip width=5,
			clip radius=8pt,
			only when rendering/.style={
			},
			]
			\node at (-0.9,0.75) {$\cp=$};
			\strand [thick] (0,0)
			to (0,1.5);
			\strand [thick] (0.5,0)
			to (0.5,1.5);
			
			\strand [thick,,-{>[scale=1,length=1.8,width=5]}] (0.5,0.75)
			to [out=0,in=-90] (0.8,0.9);
			\strand [thick] (0.8,0.9)
			to [out=90,in=0] (0.6,1.05);
			
			\strand [thick] (0,0.75)
			to [out=180,in=-90] (-0.3,0.9)
			to [out=90,in=180] (-0.1,1.05);
		\end{knot}
	\end{tikzpicture}.
\end{equation}
In the context of skein valued open Gromov-Witten theory, this map corresponds to choosing capping paths for Reeb chords of the Legendrian at infinity, hence the name. 

Observe that $\cp \circ l(\cdot,e) = P$, $\cp \circ r(e.\cdot) = \bigcirc$, and $\cp(c^m) = a C_m$ for $m \in \N$.

Encircling the annulus by a meridian above both of the boundary points
induces a map (which we denote by the same symbol as the analogous map on $\mathcal{C}$) $P: \mathcal{A} \to \mathcal{A}$, which restricts to maps $P: \mathcal{A}_n \to \mathcal{A}_n$:
\begin{equation}
	\begin{tikzpicture}
		\tikzstyle{smallnode}=[circle, inner sep=0mm, outer sep=0mm, minimum size=0.7mm, draw=black, fill=black];
		\begin{knot}[
			ignore endpoint intersections=true,
			clip width=5,
			clip radius=8pt,
			only when rendering/.style={
			},
			]
			\node at (-1.3,1.1) {$P=$};		
			\strand [thick] (0,0.5)
			to [out=90, in=-90] (0,1.7); 
			\strand [thick] (-0.4,0.5)
			to [out=90, in=-90] (-0.4,1.7); 
			
			\strand [thick, ->] (-0.44,1.52)
			to [out=190,in=90] (-0.6,1.4)
			to [out=-90,in=180] (-0.2,1.2)
			to [out=0,in=-90] (0.2,1.4)
			to [out=90,in=-10] (0.04,1.52);
			
			\node[smallnode] at (-0.8,1.1) {};
			\strand [thick] (-0.8,1.1)
			to (-0.4,1.1);
			\strand [thick] (0,1.1)
			to (0.4,1.1);
			\node[smallnode] at (0.4,1.1) {};
			\flipcrossings{1,2}
		\end{knot}
	\end{tikzpicture}.
\end{equation}
As shown above, $\left( P - \bigcirc \right)|_{\mathcal{A}_n}$ is invertible for $n \geq 1$, and it follows from the definitions that $(P - \bigcirc)|_{\mathcal{A}_0} = 0$.\\

\section{Skein valued recursion relations}\label{sec:skein valued recursion relations}

In this section, we introduce (relative) recursion relations in the skein of a solid torus. First, we define formal power series and recursion relations with values in $\mathcal{C}^+$ in Sections \ref{ssec:power series in C} and \ref{ssec:recursion in C}, and then we extend the definitions to $\mathcal{A}^+$ in Section \ref{ssec:power seires in A} and \ref{ssec:recursion in A}. The main result of this section is the existence of solutions for certain types of relative recursion relations (Theorem \ref{thm:relative recursion relation with solution}), which is stated and proved in Section \ref{ssec:existence of solutions}. After that, we discuss the behavior of relative recursion relations under applying the logarithm if it is defined (Example \ref{ex:recursion and exp/ln}), which will be play an important role in Section \ref{sec:bps partition functions}.

\subsection{Skein valued power series}\label{ssec:power series in C} Let $\mathcal{M}$ be a $\Lambda$-module. Let $\mathcal{C}^{\mathcal{M}}$ denote the $\Lambda$-module of all formal power series of the form 
\begin{equation}
\Phi(t) = \sum_{i\in\N_0} \Phi_i t^i
\end{equation}
with $\Phi_i \in \mathcal{C}_i \otimes \mathcal{M}$, where the tensor product means the tensor product over $\Lambda$. We will often omit $t$ from the notation and simply write 
\begin{equation}
\Phi = \sum_{i\in\N_0} \Phi_i.
\end{equation}

Fix an element $m \in \mathcal{M}$. We denote by $\cmm \subseteq \cm$ the subset of those $\Phi \in \cm$ with $\Phi_0 = 1 \otimes m$. 

Given a $\Lambda$-linear map $O:\mathcal{C} \to \mathcal{C}$, we will often denote the induced map $O \otimes Id_M: \mathcal{C} \otimes M \to \mathcal{C} \otimes M$ simply by $O$ as well by abuse of notation. 



Since the operator $P-\bigcirc$ has homogeneous degree zero, it induces a $\Lambda$-linear map $(P-\bigcirc) \otimes id_{\mathcal{M}}: \cm \to \cm$ which restricts to an invertible map $\cmzero \to \cmzero$.\\

Using the multiplication in $\mathcal{C}^+$, there exist natural maps
\begin{equation}
\left( \mathcal{C}_i \otimes End_\Lambda(\mathcal{M}) \right) \otimes \left( \mathcal{C}_j \otimes \mathcal{M} \right) \to \mathcal{C}_{i+j} \otimes \mathcal{M}
\end{equation}
induced by 
\begin{equation}
(\Phi_i \otimes A) \otimes (\Psi_j \otimes m) \mapsto (\Phi_i \Psi_j) \otimes A(m).
\end{equation}
These induce a map 
\begin{equation}
\begin{split}
\cendmA \otimes \cmm &\to \mathcal{C}^{\mathcal{M}}_{A(m)}\\
\left( \sum_{i \in \N_0} A_i \right) \otimes \left( \sum_{j \in \N_0} \Phi_j \right) &\mapsto \sum_{n \in \N_0} \sum_{i=0}^{n} A_i(\Phi_{n-i}),
\end{split}
\end{equation}
where we view $A_i$ as an element of $Hom_\Lambda(\mathcal{C}_{j} \otimes \mathcal{M}, \mathcal{C}_{i+j} \otimes \mathcal{M})$.\\

If $\mathcal{M}$ has the structure of a $\Lambda$-algebra, then $\cm$ inherits the structure of a $\Lambda$-algebra with multiplication induced by 
\begin{equation}
\begin{split}
\left( \mathcal{C}_i \otimes \mathcal{M} \right) \otimes \left( \mathcal{C}_j \otimes \mathcal{M} \right) &\to \mathcal{C}_{i+j} \otimes \mathcal{M}\\
(A \otimes m) \otimes (B \otimes n) &\mapsto (AB) \otimes (mn).
\end{split}
\end{equation}\\

\begin{ex}\label{ex:homflypt partition function} One important class of examples is recursively defined as $\widehat{\mathcal{C}}^l \coloneqq \mathcal{C}^{\widehat{\mathcal{C}}^{l-1}}$ with $\widehat{\mathcal{C}}^0 \coloneqq \Lambda$. Using the basis elements $W_\lambda$, we can represent any element $\Psi$ in $\widehat{\mathcal{C}}^l$ as 
\begin{equation}
	\Psi = \sum_{n_1,...,n_l\geq0} \sum_{\lambda_1 \vdash n_1}... \sum_{\lambda_l \vdash n_l} \Psi_{\lambda_1,...,\lambda_l} W_{\lambda_1} \otimes ... \otimes W_{\lambda_l} \in \widehat{\mathcal{C}}^l
\end{equation}
for some coefficients $\Psi_{\lambda_1,...,\lambda_l} \in \Lambda$. As before, we denote by $\widehat{\mathcal{C}}^l_m$ the subset of those $\Psi \in \widehat{\mathcal{C}}^l$ that satisfy $\Psi_{\emptyset,...,\emptyset} = m$.
	
	Given a framed $l$-component link $L \subseteq S^3$, we call
	\begin{equation}
		\Psi_L \coloneqq \sum_{n_1,...,n_l\geq0} \sum_{\lambda_1 \vdash n_1}... \sum_{\lambda_l \vdash n_l} L(W_{\lambda_1} \otimes ... \otimes W_{\lambda_l}) W_{\lambda_1} \otimes ... \otimes W_{\lambda_l} \in \widehat{\mathcal{C}}^l_1
	\end{equation}
	the \emph{HOMFLY-PT partition function of $L$}, where by convention we set $L(W_\emptyset,...,W_\emptyset) = 1$.\\
\end{ex}

\subsection{Recursion in $\mathcal{C}^+$}\label{ssec:recursion in C}
\begin{defn} Fix an element $m \in \mathcal{M}$. Given $\Phi \in \cmm$ and $A \in \cendmzero$, we say that \emph{$\Phi$ solves the recursion relation $P - \bigcirc = A$ (or $P - \bigcirc - A = 0$) with the initial value $m$}, or simply \emph{$\Phi$ solves $P - \bigcirc = A$}, if 
\begin{equation}
\left( (P - \bigcirc)\otimes id_{\mathcal{M}} \right)(\Phi) = A(\Phi) \in \cmzero.
\end{equation}\\
\end{defn}

Note that as $(P - \bigcirc)\otimes id_{\mathcal{M}}$ is invertible on $\cmzero$, there exists for any $A \in \cendmzero$ a unique $\Phi \in \cmm$ solving $P - \bigcirc = A$ since we can solve for the degree $n$ part $\Phi_n$ of $\Phi$ recursively (note that by assumption $A_0 = 0$):
\begin{equation}
\begin{split}
\Phi_0 &= 1 \otimes m,\\
\Phi_n &= \left((P- \bigcirc)^{-1} \otimes id_{\mathcal{M}} \right) \left(\sum_{i=1}^{n} A_i(\Phi_{n-i}) \right), \quad n \geq 1.
\end{split}
\end{equation}

Conversely, given $\Phi \in \cmm$, the values $A_n(1 \otimes m) \in \mathcal{C}_n \otimes \mathcal{M}$ for all $n \in \N$ are uniquely determined for any $A \in \cendm$ for which $\Phi$ solves $P- \bigcirc  = A$ since they satisfy
\begin{equation}
A_n(1 \otimes m) = (P - \bigcirc) \Phi_n - \sum_{i=1}^{n-1} A_i \Phi_{n-i}.
\end{equation}

In the case that $\mathcal{M} = \Lambda$ and $m = 1$, this determines $A_n \in \mathcal{C}_n \otimes End_\Lambda(\Lambda) \approx \mathcal{C}_n$ uniquely.\\

\subsection{Skein valued power series in $\mathcal{A}^+$}\label{ssec:power seires in A} Similarly to the above, we denote by $\am$ the $\Lambda$-module of formal power series of the form $\Phi = \sum_{i \in \N_0} \Phi_i$ with $\Phi_i \in \mathcal{A}_i \otimes \mathcal{M}$. For an element $m \in \mathcal{M}$, $\amm \subseteq \am$ is the subset of those $\Phi \in \am$ with $\Phi_0 = e \otimes m$.\\

The maps $l: \mathcal{C} \otimes \mathcal{A} \to \mathcal{A}$ and $r: \mathcal{A} \otimes \mathcal{C} \to \mathcal{A}$ induce  morphisms of $\Lambda$-modules (which will all be denoted by the same symbols) $l: \cm \otimes \mathcal{A}^+ \to \am$, $r: \mathcal{A}^+ \otimes \cm \to \am$, $l: \cmm \otimes \aendmA \to \mathcal{A}^{\mathcal{M}}_{A(m)}$, and $r: \aendmA \otimes \cmm \to \mathcal{A}^{\mathcal{M}}_{A(m)}$. In particular, we obtain morphisms of $\Lambda$-modules $l(\cdot,e),r(e,\cdot): \cm \to \am$. As before, their difference is denoted by $[\cdot,e] = l(\cdot,e) - r(e,\cdot)$.\\

The multiplication in $\mathcal{A}^+$ induces a map
\begin{equation}
\aendmA \otimes \amm \to \mathcal{A}^\mathcal{M}_{A(m)}
\end{equation}\\
and, if $\mathcal{M}$ is a $\Lambda$-algebra, the structure of a $\Lambda$-algebra on $\am$ and maps $l: \cm \otimes \am \to \am$ and $r: \am \otimes \cm \to \am$.\\

\subsection{Recursion in $\mathcal{A}^+$}\label{ssec:recursion in A}
\begin{defn} Fix an element $m \in \mathcal{M}$. Given $\Phi \in \cmm$ and $A \in \aendmzero$, we say that \emph{$\Phi$ solves the relative recursion relation $[\cdot,e] = r(A,\cdot)$ (resp. $[\cdot,e] = l(A,\cdot)$) with initial value $m$}, or simply \emph{$\Phi$ solves $[\cdot,e] = r(A,\cdot)$ (resp. $[\cdot,e] = l(A,\cdot)$)}, if 
	\begin{equation}
	[\Phi,e] = r(A,\Phi) \in \amzero \quad \left(\text{resp. } [\Phi,e] = l(\Phi,A) \in \amzero \right).
	\end{equation}\\
\end{defn}

Given $\Phi \in \cmm$, the values $A_n(e \otimes m) \in \mathcal{A}_n \otimes \mathcal{M}$ for $n \in \N$ are uniquely determined for any $A \in \aendm_0$ for which $\Phi$ solves $[\cdot,e] = r(A,\cdot)$ since they satisfy
\begin{equation}
A_n(e \otimes m) = r(A_n, 1 \otimes m) = [\Phi_n,e] - \sum_{i=1}^{n-1} r(A_i,\Phi_{n-i}).
\end{equation}

In the case that $\mathcal{M} = \Lambda$ and $m=1$, this determines $A_n \in \mathcal{A}_n \otimes End_\Lambda(\Lambda) \approx \mathcal{A}_n$ uniquely.\\

On the other hand, even when $\mathcal{M} = \Lambda$, not every $A \in \mathcal{A}^{\Lambda}_0$ defines a relative recursion relation that can be solved by some $\Phi \in \mathcal{C}^\Lambda_1$. Indeed, if $A_1,...,A_{n-1}$ are given, $A_n$ needs to satisfy 
\begin{equation}
A_n \in \sum_{i=1}^{n-1} r(A_i,\Phi_{n-i}) + [\mathcal{C}_n,e],
\end{equation}
but the map $[\cdot,e]: \mathcal{C}_n \to \mathcal{A}_n$ is not surjective as $\mathcal{C}_n$ has rank $p(n)$ and $\mathcal{A}_n$ has rank $\sum_{i=0}^{i=n} p(n)$, where $p(n)$ denotes the number of integer partitions of $n$.\\ 

Similar remarks apply to relative recursion relations with $l$ instead of $r$.\\

\subsection{Existence of solutions}\label{ssec:existence of solutions}
Now we will see that there are certain classes of relative recursion relations that admit a solution. For this consider the wiring

\begin{equation}\label{eq:wiring of V_n}
	\begin{tikzpicture}
	\tikzstyle{smallnode}=[circle, inner sep=0mm, outer sep=0mm, minimum size=0.7mm, draw=black, fill=black];
	\begin{knot}[
	ignore endpoint intersections=false
	clip width=5,
	clip radius=8pt,
	flip crossing=1,
	only when rendering/.style={
	}
	]
	\node at (-2.7,1.6) {$V_n = a \sum\limits_{j=1}^n$};
	\draw[draw=black] (-1.5,2) rectangle ++(1.4,0.4);
	\draw[draw=black] (0.3,2) rectangle ++(0.4,0.4);
	\node[smallnode] (left)   at (-1.8,1.6) {};
	\strand [thick] (-1.8,1.6)
	to [out=right, in=down, looseness=1.5] (-1.3,2);
	\node at (-1.1,0.3) {$n$};
	\strand [thick] (-1.1,0.6)
	to [out=up, in=down] (-1.1,2);
	\strand [thick] (-0.9,0.6)
	to [out=up, in=down] (-0.9,2);
	\strand [thick] (-0.7,0.6) node[below]{$j$}
	to [out=up,in=down, looseness=0.7] (0.7,1)
	to [out=up, in=down, looseness=0.7] (-0.7,1.4)
	to [out=up, in=down] (-0.7,2);
	\strand [thick] (-0.5,0.6)
	to [out=up, in=down] (-0.5,2);
	\strand [thick] (-0.3,0.6) node[below]{$1$}
	to [out=up, in=down] (-0.3,2);
	\strand [thick] (0.5,0.6)
	to [out=up, in=down] (0.5,2);
	\strand [thick,-{>[scale=1.7,length=2,width=2]}] (-1.3,2.4)
	to [out=up, in=down, looseness=1.5] (-1.1,2.8);
	\strand [thick,-{>[scale=1.7,length=2,width=2]}] (-1.1,2.4)
	to [out=up, in=down, looseness=1.5] (-0.9,2.8);
	\strand [thick,-{>[scale=1.7,length=2,width=2]}] (-0.9,2.4)
	to [out=up, in=down, looseness=1.5] (-0.7,2.8);
	\strand [thick,-{>[scale=1.7,length=2,width=2]}] (-0.7,2.4)
	to [out=up, in=down, looseness=1.5] (-0.5,2.8);
	\strand [thick,-{>[scale=1.7,length=2,width=2]}] (-0.5,2.4)
	to [out=up, in=down, looseness=1.5] (-0.3,2.8);
	\strand [thick,-{>[scale=1.7,length=2,width=2]}] (0.5,2.4)
	to [out=up, in=down, looseness=1.5] (0.5,2.8);
	\node[smallnode] (right)   at (0.9,1.6) {};
	\strand [thick] (-0.3,2.4)
	to [out=up, in=up, looseness=1.5] (0.1,2.4)
	to [out=down, in = up] (0.1,1.8)
	to [out=down, in = left, looseness = 0.8] (0.9,1.6);
	
	\node at (1.6,1.6) {$- \sum\limits_{j=1}^n$};
	\draw[draw=black] (2.5,2) rectangle ++(1.4,0.4);
	\draw[draw=black] (4.3,2) rectangle ++(0.4,0.4);
	\node[smallnode] (left)   at (2.2,1.6) {};
	\strand [thick] (2.2,1.6)
	to [out=right, in=left] (2.8,1.6)
	to [out=right, in=down, looseness=1.5] (3.3,2);
	\node at (2.9,0.3) {$n$};
	\strand [thick] (2.9,0.6)
	to [out=up, in=down] (2.9,2);
	\strand [thick] (3.1,0.6)
	to [out=up, in=down] (3.1,2);
	\strand [thick] (3.7,2.4)
	to [out=up, in=right, looseness=1.5] (3.45,2.6)
	to [out=left, in=right] (2.7,2.6)
	to [out=left, in=up,looseness=1.5] (2.4,2.3)
	to [out=down, in=up] (2.4,2)
	to [out=down, in=down, looseness=1.5] (2.7,2);
	\strand [thick] (3.5,0.6)
	to [out=up, in=down] (3.5,2);
	\strand [thick] (3.7,0.6) node[below]{$1$}
	to [out=up, in=down] (3.7,2);
	\strand [thick] (4.5,0.6)
	to [out=up, in=down] (4.5,2);
	\strand [thick,-{>[scale=1.7,length=2,width=2]}] (2.7,2.4)
	to [out=up, in=down, looseness=1.5] (2.9,2.8);
	\strand [thick,-{>[scale=1.7,length=2,width=2]}] (2.9,2.4)
	to [out=up, in=down, looseness=1.5] (3.1,2.8);
	\strand [thick,-{>[scale=1.7,length=2,width=2]}] (3.1,2.4)
	to [out=up, in=down, looseness=1.5] (3.3,2.8);
	\strand [thick,-{>[scale=1.7,length=2,width=2]}] (3.3,2.4)
	to [out=up, in=down, looseness=1.5] (3.5,2.8);
	\strand [thick,-{>[scale=1.7,length=2,width=2]}] (3.5,2.4)
	to [out=up, in=down, looseness=1.5] (3.7,2.8);
	\strand [thick,-{>[scale=1.7,length=2,width=2]}] (4.5,2.4)
	to [out=up, in=down, looseness=1.5] (4.5,2.8);
	\node[smallnode] (right)   at (4.9,1.6) {};
	\strand [thick] (3.3,0.6) node[below]{$j$}
	to [out=up,in=down] (3.3,1.2)
	to [out=up, in=left, looseness=1.5] (3.8,1.6)
	to [out=right, in=left] (4.9,1.6);
	
	\flipcrossings{2,4,6,8,9,10,11,12,13,14,15,16,17}
	\end{knot}
	\end{tikzpicture},
\end{equation}

and recall that

\begin{equation}\label{eq:wiring W_n2}
\begin{tikzpicture}
\tikzstyle{smallnode}=[circle, inner sep=0mm, outer sep=0mm, minimum size=0.7mm, draw=black, fill=black];
\begin{knot}[
ignore endpoint intersections=false
clip width=5,
clip radius=8pt,
flip crossing=1,
only when rendering/.style={
}
]
\node at (-2.5,1.6) {$W_n =$};
\draw[draw=black] (-1.5,2) rectangle ++(1.4,0.4);
\node[smallnode] (left)   at (-1.8,1.6) {};
\strand [thick] (-1.8,1.6)
to [out=right, in=down, looseness=1.5] (-1.3,2);
\node at (-1.1,0.3) {$n$};
\strand [thick] (-1.1,0.6)
to [out=up, in=down] (-1.1,2);
\strand [thick] (-0.9,0.6)
to [out=up, in=down] (-0.9,2);
\strand [thick] (-0.7,0.6)
to [out=up, in=down] (-0.7,2);
\strand [thick] (-0.5,0.6)
to [out=up, in=down] (-0.5,2);
\strand [thick] (-0.3,0.6) node[below]{$1$}
to [out=up, in=down] (-0.3,2);
\strand [thick,-{>[scale=1.7,length=2,width=2]}] (-1.3,2.4)
to [out=up, in=down, looseness=1.5] (-1.1,2.8);
\strand [thick,-{>[scale=1.7,length=2,width=2]}] (-1.1,2.4)
to [out=up, in=down, looseness=1.5] (-0.9,2.8);
\strand [thick,-{>[scale=1.7,length=2,width=2]}] (-0.9,2.4)
to [out=up, in=down, looseness=1.5] (-0.7,2.8);
\strand [thick,-{>[scale=1.7,length=2,width=2]}] (-0.7,2.4)
to [out=up, in=down, looseness=1.5] (-0.5,2.8);
\strand [thick,-{>[scale=1.7,length=2,width=2]}] (-0.5,2.4)
to [out=up, in=down, looseness=1.5] (-0.3,2.8);
\node[smallnode] (right)   at (0.9,1.6) {};
\strand [thick] (-0.3,2.4)
to [out=up, in=up, looseness=1.5] (0.1,2.4)
to [out=down, in = up] (0.1,1.8)
to [out=down, in = left, looseness = 0.8] (0.9,1.6);

\end{knot}
\end{tikzpicture},
\end{equation}
which induce maps $V_n: R^{n+1}_{n+1} \to \Hom_\Lambda(\mathcal{C}_n,\mathcal{A}_n)$ and $W_n: R^{n+1}_{n+1} \to \mathcal{A}_n$, where in each summand in the definition of $V_n$, the wide rectangle is diffeomorphic to $[0,1]^2$ with $2(n+1)$ boundary points, which are oriented to point into the rectangle at the bottom and outwards at the top, and the square with the vertical lines ending at it represents an embedded solid torus.

\begin{thm}\label{thm:relative recursion relation with solution}
	For every $n \geq 1$ fix a finite set $I_n$. Let $a_\kappa \in R^{n+1}_{n+1}$ and $A_\kappa \in End_\Lambda(\mathcal{M}), \kappa \in I_n,$ be such that $[A_\kappa,A_\chi] = 0$ for all $i,j \in \N$, $\kappa \in I_i$, and $\chi \in I_j$, and
	\begin{equation}\label{eq:V_n(a_kappa)=0}
	\sum\limits_{\kappa \in I_n} \left(V_n(a_\kappa) \otimes A_\kappa \right) = 0 \in \Hom_\Lambda(\mathcal{C}_n \otimes \mathcal{M}, \mathcal{A}_n \otimes \mathcal{M})
	\end{equation}
	for all $n \in \N$.
	
	Then the relative recursion relation $[\cdot,e] = r(A,\cdot)$, where
	\begin{equation}
	A \coloneqq \sum\limits_{n=1}^\infty \sum\limits_{\kappa \in I_n} W_n(a_\kappa) \otimes A_\kappa,
	\end{equation}
	admits a solution in $\cmm$ for any $m \in \mathcal{M}$.
\end{thm}

\bproof Let $\Phi = \sum_{n \in \N_0} \Phi_n \in \cmm$ be the solution of $P - \bigcirc = \left(\cp \otimes id_{\mathcal{M}}\right)(A) \in \cendmzero$ with initial value $m$. We will show that $\Phi$ also solves the relative recursion relation $[\cdot,e] = r(A,\cdot)$ by using induction on $n$.\\

For $n=1$, we may assume without loss of generality that $I_1 = \{1,2\}$ and $a_1 = \kappa_1 1_2, a_2 = \kappa_2 \sigma_1 \in R^2_2$ for some $\kappa_1,\kappa_2 \in \Lambda,$ which are mapped by $W_1$ to $\kappa_1 c$ and $\kappa_2 l(W_{(1)},e)$, where $1_2$ denotes the trivial braid and $\sigma_1$ denotes the braid with one positive crossing. Using the skein relations, one obtains that $V_1(1_2) = 0$ and $V_1(\sigma_1)(1) = a l(W_{(1)},e) - (\bigcirc + az) c$. It follows from (\ref{eq:V_n(a_kappa)=0}) and the fact that $l(W_{(1)},e)$ and $c$ are linearly independent over $\Lambda$ that $\kappa_2 l(W_{(1)},e) \otimes A_2 = 0$, which implies that 
\begin{equation}
\sum\limits_{\kappa \in I_1} W_n(a_\kappa) \otimes A_\kappa = \kappa_1 c \otimes A_1,
\end{equation}
which is mapped by $\cp \otimes id_{\mathcal{M}}$ to $a \kappa_1 W_{(1)} \otimes A_1$. Thus, we have that 
\begin{equation}
(P - \bigcirc)(\Phi_1) = a \kappa_1 W_{(1)} \otimes A_1(m),
\end{equation}
which implies that $\Phi_1 = z^{-1} \kappa_1 W_{(1)} \otimes A_1(m)$. It follows that 
\begin{equation}
[\Phi_1,e] = \kappa_1 c \otimes A_1(m) = r\left( \sum\limits_{\kappa \in I_1} W_n(a_\kappa) \otimes A_\kappa, 1 \otimes m \right),
\end{equation}
as was claimed.

Now assume that we know that for some $n \in \N_{\geq 2}$, $\Phi$ satisfies 
\begin{equation}
[\Phi_k,e] = \sum\limits_{i=1}^{k} \sum\limits_{\kappa \in I_i} r\left(W_i(a_\kappa)\otimes A_\kappa,\Phi_{k-i}\right)
\end{equation}
for all $k \in \{1,...,n-1\}$. We will show that 
\begin{equation}\label{eq:p-unknot of thm}
(P - \bigcirc)\left([\Phi_n,e]\right) = (P-\bigcirc)\left(\sum\limits_{i=1}^{n} \sum\limits_{\kappa \in I_i} r\left(W_i(a_\kappa)\otimes A_\kappa,\Phi_{n-i}\right) \right),
\end{equation}
which finishes the induction step by injectivity of $P - \bigcirc$ on $\mathcal{A}_n$. Using that $(P - \bigcirc) \left( [\cdot,e] \right) = [(P-\bigcirc)(\cdot),e]$, we compute
\begin{equation}\label{eq:p-unknot of relative recursion}
\begin{split}
(P - \bigcirc)\left([\Phi_n,e]\right) &= [(P-\bigcirc)(\Phi_n),e] \\
&= \left[ \sum\limits_{i=1}^{n} \sum\limits_{\kappa \in I_i} \left(\overline{W_i(a_\kappa)} \otimes A_\kappa\right) \left(\Phi_{n-i}\right),e\right] \\
&= \sum\limits_{i=1}^{n} \sum\limits_{\kappa \in I_i} (1 \otimes A_\kappa) \left( l\left(\overline{W_i(a_\kappa)},[\Phi_{n-i},e]\right) + r\left( \left[\overline{W_i(a_\kappa)},e \right], \Phi_{n-i} \right) \right)
\end{split}
\end{equation}
where we used the notation $\overline{W_i(a_\kappa)} \coloneqq \cp(W_i(a_\kappa))$.

Using the induction hypothesis, it follows that 
\begin{equation}\label{eq:p-unknot of relative recursion1}
\begin{gathered}
\sum\limits_{i=1}^{n} \sum\limits_{\kappa \in I_i} (1 \otimes A_\kappa) \left( l\left(\overline{W_i(a_\kappa)},[\Phi_{n-i},e]\right) \right)\\ 
= \sum\limits_{i=1}^{n} \sum\limits_{\kappa \in I_i} (1 \otimes A_\kappa)\left( l\left(\overline{W_i(a_\kappa)}, \sum_{k=1}^{n-i} \sum_{\chi \in I_k} r\left(W_k(a_\chi) \otimes A_\chi,\Phi_{n-i-k}\right)\right) \right)  \\
= \underset{i+k \leq n}{\sum\limits_{i,k\geq1}} \sum\limits_{\kappa \in I_i} \sum\limits_{\chi \in I_k} (1 \otimes A_\kappa \circ A_\chi)\left( l\left(\overline{W_i(a_\kappa)},r\left(W_k(a_\chi),\Phi_{n-i-k}\right)\right) \right).
\end{gathered}
\end{equation}

For $\kappa \in I_k, k \in \N,$ denote by $b_\kappa \in \mathcal{C}_k$ the image of $a_\kappa$ under the map induced by the wiring
\begin{equation}\label{eq:wiring for b_lambda}
\begin{tikzpicture}
\tikzstyle{smallnode}=[circle, inner sep=0mm, outer sep=0mm, minimum size=0.7mm, draw=black, fill=black];
\begin{knot}[
ignore endpoint intersections=false
clip width=5,
clip radius=8pt,
flip crossing=1,
only when rendering/.style={
}
]

\node at (1.8,1.6) {$\sum\limits_{j=1}^k$};
\draw[draw=black] (2.5,2) rectangle ++(1.4,0.4);
\node[smallnode] (left)   at (2.2,1.6) {};
\strand [thick] (2.2,1.6)
to [out=right, in=left] (2.8,1.6)
to [out=right, in=down, looseness=1.5] (3.3,2);
\node at (2.9,0.3) {$k$};
\strand [thick] (2.9,0.6)
to [out=up, in=down] (2.9,2);
\strand [thick] (3.1,0.6)
to [out=up, in=down] (3.1,2);
\strand [thick] (3.7,2.4)
to [out=up, in=right, looseness=1.5] (3.45,2.6)
to [out=left, in=right] (2.7,2.6)
to [out=left, in=up,looseness=1.5] (2.4,2.3)
to [out=down, in=up] (2.4,2)
to [out=down, in=down, looseness=1.5] (2.7,2);
\strand [thick] (3.5,0.6)
to [out=up, in=down] (3.5,2);
\strand [thick] (3.7,0.6) node[below]{$1$}
to [out=up, in=down] (3.7,2);
\strand [thick,-{>[scale=1.7,length=2,width=2]}] (2.7,2.4)
to [out=up, in=down, looseness=1.5] (2.9,2.8);
\strand [thick,-{>[scale=1.7,length=2,width=2]}] (2.9,2.4)
to [out=up, in=down, looseness=1.5] (3.1,2.8);
\strand [thick,-{>[scale=1.7,length=2,width=2]}] (3.1,2.4)
to [out=up, in=down, looseness=1.5] (3.3,2.8);
\strand [thick,-{>[scale=1.7,length=2,width=2]}] (3.3,2.4)
to [out=up, in=down, looseness=1.5] (3.5,2.8);
\strand [thick,-{>[scale=1.7,length=2,width=2]}] (3.5,2.4)
to [out=up, in=down, looseness=1.5] (3.7,2.8);
\node[smallnode] (right)   at (4.2,1.6) {};
\strand [thick] (3.3,0.6) node[below]{$j$}
to [out=up,in=down] (3.3,1.2)
to [out=up, in=left, looseness=1.5] (3.8,1.6)
to [out=right, in=left] (4.2,1.6);

\flipcrossings{1,2,3,4,5,6,7,8,9}
\end{knot}
\end{tikzpicture}.
\end{equation}

Using the skein relations, one computes that
\begin{equation}\label{eq:p-unknot of relative recursion2}
	\left[\overline{W_i(a_\kappa)},e \right] = z b_\kappa
\end{equation}
for $i \in \N, \kappa \in I_i$ by successively passing the strand of $e$ through the parallel vertical strands in the diagram of $\overline{W_i(a_\kappa)}$.

It follows that
\begin{equation}\label{eq:p-unknot of relative recursion1inin}
l\left(\overline{W_i(a_\kappa)},W_k(a_\chi) \right) = l\left(\overline{W_i(a_\kappa)},e \right) \cdot W_k(a_\chi) = r\left(W_k(a_\chi),\overline{W_i(a_\kappa)} \right) + z b_\kappa \cdot W_k(a_\chi)
\end{equation}
for $i,k \in \N, \kappa \in I_i, \chi \in I_k$ where $``\cdot"$ denotes the product of the $\Lambda$-algebra $\mathcal{A}$. Using this, we compute further that
\begin{equation}\label{eq:p-unknot of relative recursion1in}
\begin{split}
l\left(\overline{W_i(a_\kappa)},r\left(W_k(a_\chi),\Phi_{n-i-k}\right)\right) &= r\left(l\left(\overline{W_i(a_\kappa)},W_k(a_\chi) \right), \Phi_{n-i-k} \right) \\
&= r\left(W_k(a_\chi), \overline{W_i(a_\kappa)} \cdot \Phi_{n-i-k} \right) + z b_\kappa \cdot r\left( W_k(a_\chi), \Phi_{n-i-k} \right),
\end{split}
\end{equation}
where in the first summand of the second row, $``\cdot"$ denotes the product of $\mathcal{C}$.

We can write $A_\kappa \circ A_\chi = [A_\kappa,A_\chi] + A_\chi \circ A_\kappa$, and after applying $A_\chi \circ A_\kappa$ and summing, the first summand of the second row in (\ref{eq:p-unknot of relative recursion1in}) becomes
\begin{equation}\label{eq:p-unknot of relative recursion1in2}
\begin{gathered}
\underset{i+k \leq n}{\sum\limits_{i,k\geq1}} \sum\limits_{\kappa \in I_i} \sum\limits_{\chi \in I_k} (1 \otimes A_\chi \circ A_\kappa)\left( r\left(W_k(a_\chi), \overline{W_i(a_\kappa)} \cdot \Phi_{n-i-k} \right) \right)\\
= \sum\limits_{k = 1}^n \sum\limits_{\chi \in I_k} (1 \otimes A_\chi)\left( r\left(W_k(a_\chi), \sum_{i=1}^{n-k} \sum\limits_{\kappa \in I_i} \left( \overline{W_i(a_\kappa)} \otimes A_\kappa \right) \left(\Phi_{n-i-k} \right) \right) \right)\\
= \sum\limits_{k = 1}^n \sum\limits_{\chi \in I_k} (1 \otimes A_\chi)\left( r\left(W_k(a_\chi), \left( P - \bigcirc \right) \left(\Phi_{n-k} \right) \right) \right)\\
= \sum\limits_{i = 1}^n \sum\limits_{\kappa \in I_i} (1 \otimes A_\kappa)\left( r\left(W_i(a_\kappa), \left( P - \bigcirc \right) \left(\Phi_{n-i} \right) \right) \right),
\end{gathered}
\end{equation}
where for the second equality we used the definition of $\Phi$ and in the last step we simply renamed the summation variables.

For $\kappa \in I_i, i \in \N$, let $c_\kappa \in \Hom_\Lambda(\mathcal{C}_{n-i},\mathcal{A}_{n})$ be the image of $a_\kappa$ under the map induced by the wiring
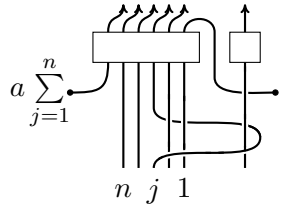
\begin{equation}\label{eq:wiring for c_lambda}
\begin{tikzpicture}
\tikzstyle{smallnode}=[circle, inner sep=0mm, outer sep=0mm, minimum size=0.7mm, draw=black, fill=black];
\begin{knot}[
ignore endpoint intersections=false
clip width=5,
clip radius=8pt,
flip crossing=1,
only when rendering/.style={
}
]
\node at (-2.2,1.6) {$a \sum\limits_{j=1}^n$};
\draw[draw=black] (-1.5,2) rectangle ++(1.4,0.4);
\draw[draw=black] (0.3,2) rectangle ++(0.4,0.4);
\node[smallnode] (left)   at (-1.8,1.6) {};
\strand [thick] (-1.8,1.6)
to [out=right, in=down, looseness=1.5] (-1.3,2);
\node at (-1.1,0.3) {$n$};
\strand [thick] (-1.1,0.6)
to [out=up, in=down] (-1.1,2);
\strand [thick] (-0.9,0.6)
to [out=up, in=down] (-0.9,2);
\strand [thick] (-0.7,0.6) node[below]{$j$}
to [out=up,in=down, looseness=0.7] (0.7,1)
to [out=up, in=down, looseness=0.7] (-0.7,1.4)
to [out=up, in=down] (-0.7,2);
\strand [thick] (-0.5,0.6)
to [out=up, in=down] (-0.5,2);
\strand [thick] (-0.3,0.6) node[below]{$1$}
to [out=up, in=down] (-0.3,2);
\strand [thick] (0.5,0.6)
to [out=up, in=down] (0.5,2);
\strand [thick,-{>[scale=1.7,length=2,width=2]}] (-1.3,2.4)
to [out=up, in=down, looseness=1.5] (-1.1,2.8);
\strand [thick,-{>[scale=1.7,length=2,width=2]}] (-1.1,2.4)
to [out=up, in=down, looseness=1.5] (-0.9,2.8);
\strand [thick,-{>[scale=1.7,length=2,width=2]}] (-0.9,2.4)
to [out=up, in=down, looseness=1.5] (-0.7,2.8);
\strand [thick,-{>[scale=1.7,length=2,width=2]}] (-0.7,2.4)
to [out=up, in=down, looseness=1.5] (-0.5,2.8);
\strand [thick,-{>[scale=1.7,length=2,width=2]}] (-0.5,2.4)
to [out=up, in=down, looseness=1.5] (-0.3,2.8);
\strand [thick,-{>[scale=1.7,length=2,width=2]}] (0.5,2.4)
to [out=up, in=down, looseness=1.5] (0.5,2.8);
\node[smallnode] (right)   at (0.9,1.6) {};
\strand [thick] (-0.3,2.4)
to [out=up, in=up, looseness=1.5] (0.1,2.4)
to [out=down, in = up] (0.1,1.8)
to [out=down, in = left, looseness = 0.8] (0.9,1.6);

\flipcrossings{2,4,6}
\end{knot}
\end{tikzpicture}.
\end{equation}

Using the skein relations, it follows that
\begin{equation}
(P-\bigcirc)\left(r\left( W_i(a_\kappa),\Phi_{n-i} \right) \right) - r\left( W_i(a_\kappa),(P-\bigcirc)(\Phi_{n-i}) \right) = z c_\kappa(\Phi_{n-i})
\end{equation}
by passing the meridian curve successively through the the vertical parallel strands in the diagram of $W_i(a_\kappa)$.

With this, we can rewrite (\ref{eq:p-unknot of relative recursion1in2}) as
\begin{equation}\label{eq:p-unknot of relative recursion1in2.2}
\begin{gathered}
\sum\limits_{i = 1}^n \sum\limits_{\kappa \in I_k} (1 \otimes A_\kappa)\left( r\left(W_i(a_\kappa), \left( P - \bigcirc \right) \left(\Phi_{n-i} \right) \right) \right)\\
= (P-\bigcirc)\left(\sum\limits_{i=1}^{n} \sum\limits_{\kappa \in I_i} r\left(W_i(a_\kappa)\otimes A_\kappa,\Phi_{n-i}\right) \right) - z \sum\limits_{i = 1}^n \sum\limits_{\kappa \in I_k} \left( c_\kappa \otimes A_\kappa \right) (\Phi_{n-i}).
\end{gathered}
\end{equation}

All in all, we obtain that 
\begin{equation}\label{eq:p-unknot of relative recursion endresult}
\begin{gathered}
(P - \bigcirc)\left([\Phi_n,e]\right) \overset{(\ref{eq:p-unknot of relative recursion}),(\ref{eq:p-unknot of relative recursion1})}{=} \underset{i+k \leq n}{\sum\limits_{i,k\geq1}} \sum\limits_{\kappa \in I_i} \sum\limits_{\chi \in I_k} (1 \otimes A_\kappa \circ A_\chi)\left( l\left(\overline{W_i(a_\kappa)},r\left(W_k(a_\chi),\Phi_{n-i-k}\right)\right) \right)\\ + \sum\limits_{i=1}^{n} \sum\limits_{\kappa \in I_i} (1 \otimes A_\kappa) \left(r\left( \left[\overline{W_i(a_\kappa)},e \right], \Phi_{n-i} \right) \right)\\
\overset{(\ref{eq:p-unknot of relative recursion1in}),(\ref{eq:p-unknot of relative recursion2})}{=} \underset{i+k \leq n}{\sum\limits_{i,k\geq1}} \sum\limits_{\kappa \in I_i} \sum\limits_{\chi \in I_k} (1 \otimes A_\kappa \circ A_\chi)\left( r\left(W_k(a_\chi), \overline{W_i(a_\kappa)} \cdot \Phi_{n-i-k} \right) + z b_\kappa \cdot r\left( W_k(a_\chi), \Phi_{n-i-k} \right) \right)\\
+ z \sum\limits_{i=1}^{n} \sum\limits_{\kappa \in I_i} (1 \otimes A_\kappa) \left(r\left( b_\kappa, \Phi_{n-i} \right) \right)\\
\overset{(\ref{eq:p-unknot of relative recursion1in2}),(\ref{eq:p-unknot of relative recursion1in2.2})}{=} (P-\bigcirc)\left(\sum\limits_{i=1}^{n} \sum\limits_{\kappa \in I_i} r\left(W_i(a_\kappa)\otimes A_\kappa,\Phi_{n-i}\right) \right) \\
- z \sum\limits_{i=1}^{n} \sum\limits_{\kappa \in I_i} \left(b_\kappa \otimes A_\kappa \right) \left([\Phi_{n-i},e] - r\left(\sum_{k=1}^{n-i} \sum_{\chi \in I_k} W_k(a_\chi) \otimes A_\chi, \Phi_{n-i-k} \right) \right) \\
+ z \sum\limits_{i=1}^{n} \sum\limits_{\kappa \in I_i} (1 \otimes A_\kappa) \left(l(\Phi_{n-i},b_\kappa) - c_\kappa(\Phi_{n-i}) \right)\\
+ \underset{i+k \leq n}{\sum\limits_{i,k\geq1}} \sum\limits_{\kappa \in I_i} \sum\limits_{\chi \in I_k} (1 \otimes [A_\kappa,A_\chi])\left( r\left(W_k(a_\chi), \overline{W_i(a_\kappa)} \cdot \Phi_{n-i-k} \right) \right),
\end{gathered}
\end{equation}
where in the last step we also rewrote 
\begin{equation} (1 \otimes A_\kappa) \left(r\left( b_\kappa, \Phi_{n-i} \right) \right) = - \left( b_\kappa \otimes A_\kappa \right) \left( [\Phi_{n-i},e] \right) + (1 \otimes A_\kappa) \left(l(\Phi_{n-i},b_\kappa) \right).
\end{equation}

The second term on the right hand side in (\ref{eq:p-unknot of relative recursion endresult}) vanishes by the induction hypothesis, the third term vanishes since 
\begin{equation}\sum\limits_{\kappa \in I_i} (1 \otimes A_\kappa)\left(l(\Phi_{n-i},b_\kappa) - c_\kappa(\Phi_{n-i}) \right) = -\sum\limits_{\kappa \in I_i} \left(V_i(a_\kappa) \otimes A_\kappa \right)(\Phi_{n-i}) = 0
\end{equation}
by definition of $b_\kappa$, $c_\kappa$, and $V_i$, and the last term vanishes since $[A_\kappa,A_\chi] = 0$. This finishes the proof.
\eproof\\

\begin{rmk}
Note that in the derivation of (\ref{eq:p-unknot of relative recursion endresult}) we did not use the assumptions on $a_\kappa$ and $A_\kappa$ in Theorem \ref{thm:relative recursion relation with solution}. Thus, if we are given any $a_\kappa$ and $A_\kappa$ as in Theorem \ref{thm:relative recursion relation with solution} without assuming that $[A_\kappa,A_\chi] = 0$ or $\sum V_n(a_\kappa) \otimes A_\kappa = 0$, and we assume that  
\begin{equation}
[\Phi_k,e] = \sum\limits_{i=1}^{k} \sum\limits_{\kappa \in I_i} r\left(W_i(a_\kappa)\otimes A_\kappa,\Phi_{k-i}\right)
\end{equation}
holds for all $k < n$ for some $n \in \N$, where $\Phi$ is as above the solution to the recursion relation $P - \bigcirc = \left(\cp \otimes id_M\right)(A)$ with initial condition $m \in \mathcal{M}$, then 
\begin{equation}\label{eq:relative recursion in deg n}
[\Phi_n,e] = \sum\limits_{i=1}^{n} \sum\limits_{\kappa \in I_i} r\left(W_i(a_\kappa)\otimes A_\kappa,\Phi_{n-i}\right)
\end{equation}
holds if and only if 
\begin{equation}
\begin{gathered}
[(P - \bigcirc)(\Phi_n),e] = (P-\bigcirc)\left(\sum\limits_{i=1}^{n} \sum\limits_{\kappa \in I_i} r\left(W_i(a_\kappa)\otimes A_\kappa,\Phi_{n-i}\right) \right)\\
- z \sum\limits_{i=1}^{n} \sum\limits_{\kappa \in I_i} (1 \otimes A_\kappa) \left(V_i(a_\kappa)(\Phi_{n-i}) \right)\\
+ \underset{i+k \leq n}{\sum\limits_{i,k\geq1}} \sum\limits_{\kappa \in I_i} \sum\limits_{\chi \in I_k} (1 \otimes [A_\kappa,A_\chi])\left( r\left(W_k(a_\chi), \overline{W_i(a_\kappa)} \cdot \Phi_{n-i-k} \right) \right).
\end{gathered}
\end{equation}
This gives a relation in terms of $a_\chi, A_\chi$, and $\Phi_k$, $k < n,\chi \in I_k$, that $a_\kappa$ and $A_\kappa, \kappa \in I_n,$ need to satisfy in order for (\ref{eq:relative recursion in deg n}) to hold.\\
\end{rmk}

The next result shows that $\ker V_n$ is non-empty and gives explicit elements in the kernel.

\begin{prop}\label{prop:kernel of V_n}
	For all $n \geq 1$, $1_{n+1} \in \ker V_n$.
\end{prop}


\bproof We show that, in fact, the relation

\begin{equation}\label{eq:V_n(1_{n+1})=0}
	\begin{tikzpicture}
	\tikzstyle{smallnode}=[circle, inner sep=0mm, outer sep=0mm, minimum size=0.7mm, draw=black, fill=black];
	\begin{knot}[
	end tolerance=.1pt
	clip width=5,
	clip radius=8pt,
	flip crossing=1,
	only when rendering/.style={
	}
	]
	\node at (-2.1,1.6) {$a$};
	\draw[draw=black] (0.3,2) rectangle ++(0.4,0.4);
	\node[smallnode] (left)   at (-1.8,1.6) {};
	\strand [thick,-{>[scale=1.7,length=2,width=2]}] (-1.8,1.6)
	to [out=right, in=down, looseness=1.5] (-1.3,2)
	to [out=up, in=down] (-1.3,2.4)
	to [out=up, in=down, looseness=1.5] (-1.1,2.8);
	\node at (-1.1,0.3) {$n$};
	\strand [thick,-{>[scale=1.7,length=2,width=2]}] (-1.1,0.6)
	to [out=up, in=down] (-1.1,2.4)
	to [out=up, in=down, looseness=1.5] (-0.9,2.8);
	\strand [thick,-{>[scale=1.7,length=2,width=2]}] (-0.9,0.6)
	to [out=up, in=down] (-0.9,2.4)
	to [out=up, in=down, looseness=1.5] (-0.7,2.8);
	\strand [thick,-{>[scale=1.7,length=2,width=2]}] (-0.7,0.6) node[below]{$j$}
	to [out=up,in=down, looseness=0.7] (0.7,1)
	to [out=up, in=down, looseness=0.7] (-0.7,1.4)
	to [out=up, in=down] (-0.7,2.4)
	to [out=up, in=down, looseness=1.5] (-0.5,2.8);
	\strand [thick,-{>[scale=1.7,length=2,width=2]}] (-0.5,0.6)
	to [out=up, in=down] (-0.5,2.4)
	to [out=up, in=down, looseness=1.5] (-0.3,2.8);
	\node[smallnode] (right)   at (0.9,1.6) {};
	\strand [thick] (-0.3,0.6) node[below]{$1$}
	to [out=up, in=down] (-0.3,2.4)
	to [out=up, in=up, looseness=1.5] (0.1,2.4)
	to [out=down, in = up] (0.1,1.8)
	to [out=down, in = left, looseness = 0.8] (0.9,1.6);
	\strand [thick] (0.5,0.6)
	to [out=up, in=down] (0.5,2);
	\strand [thick,-{>[scale=1.7,length=2,width=2]}] (0.5,2.4)
	to [out=up, in=down, looseness=1.5] (0.5,2.8);
	
	\node at (1.6,1.6) {$=$};
	\draw[draw=black] (4.3,2) rectangle ++(0.4,0.4);
	\node[smallnode] (left)   at (2.2,1.6) {};
	\strand [thick,-{>[scale=1.7,length=2,width=2]}] (2.2,1.6)
	to [out=right, in=left] (2.8,1.6)
	to [out=right, in=down, looseness=1.5] (3.3,2)
	to [out=up, in=down] (3.3,2.4)
	to [out=up, in=down, looseness=1.5] (3.5,2.8);
	\node at (2.9,0.3) {$n$};
	\strand [thick,-{>[scale=1.7,length=2,width=2]}] (2.9,0.6)
	to [out=up, in=down] (2.9,2.4)
	to [out=up, in=down, looseness=1.5] (3.1,2.8);
	\strand [thick,-{>[scale=1.7,length=2,width=2]}] (3.1,0.6)
	to [out=up, in=down] (3.1,2.4)
	to [out=up, in=down, looseness=1.5] (3.3,2.8);
	\node[smallnode] (right)   at (4.9,1.6) {};
	\strand [thick] (3.3,0.6) node[below]{$j$}
	to [out=up,in=down] (3.3,1.2)
	to [out=up, in=left, looseness=1.5] (3.8,1.6)
	to [out=right, in=left] (4.9,1.6);
	\strand [thick,-{>[scale=1.7,length=2,width=2]}] (3.5,0.6)
	to [out=up, in=down] (3.5,2.4)
	to [out=up, in=down, looseness=1.5] (3.7,2.8);
	\strand [thick] (3.7,0.6) node[below]{$1$}
	to [out=up, in=down] (3.7,2.4)
	to [out=up, in=right, looseness=1.5] (3.45,2.6)
	to [out=left, in=right] (2.7,2.6)
	to [out=left, in=up,looseness=1.5] (2.4,2.3);
	\strand [thick,-{>[scale=1.7,length=2,width=2]}] (2.4,2.3)
	to [out=down, in=up] (2.4,2)
	to [out=down, in=down, looseness=1.5] (2.7,2)
	to [out=up, in=down] (2.7,2.4)
	to [out=up, in=down, looseness=1.5] (2.9,2.8);
	\strand [thick] (4.5,0.6)
	to [out=up, in=down] (4.5,2);
	\strand [thick,-{>[scale=1.7,length=2,width=2]}] (4.5,2.4)
	to [out=up, in=down, looseness=1.5] (4.5,2.8);
	\flipcrossings{1,2,4,6,7,8,9,17}
	\end{knot}
	\end{tikzpicture}
\end{equation}
holds for every $j = 1,...,n$ by using induction on $j$.

In the case that $j=1$, the claim follows by a framed isotopy and one of the Skein relations. Now assume that we have proven the claim for $j-1$ for some $j \geq 2$. The left hand side of (\ref{eq:V_n(1_{n+1})=0}) can be written as follows:

\begin{equation}\label{eq:lhs of V_n(1_{n+1})=0}
\begin{tikzpicture}
\tikzstyle{smallnode}=[circle, inner sep=0mm, outer sep=0mm, minimum size=0.7mm, draw=black, fill=black];
\begin{knot}[
end tolerance=.1pt
clip width=5,
clip radius=8pt,
only when rendering/.style={
}
]
\node at (-3.1,1.6) {$a^{-1} LHS_j \coloneqq$};
\draw[draw=black] (0.3,2) rectangle ++(0.4,0.4);
\node[smallnode] (left)   at (-1.8,1.6) {};
\strand [thick,-{>[scale=1.7,length=2,width=2]}] (-1.8,1.6)
to [out=right, in=down, looseness=1.5] (-1.3,2)
to [out=up, in=down] (-1.3,2.4)
to [out=up, in=down, looseness=1.5] (-1.1,2.8);
\node at (-1.1,0.3) {$n$};
\strand [thick,-{>[scale=1.7,length=2,width=2]}] (-1.1,0.6)
to [out=up, in=down] (-1.1,2.4)
to [out=up, in=down, looseness=1.5] (-0.9,2.8);
\strand [thick,-{>[scale=1.7,length=2,width=2]}] (-0.9,0.6)
to [out=up, in=down] (-0.9,2.4)
to [out=up, in=down, looseness=1.5] (-0.7,2.8);
\strand [thick,-{>[scale=1.7,length=2,width=2]}] (-0.7,0.6) node[below]{$j$}
to [out=up,in=down, looseness=0.7] (0.7,1)
to [out=up, in=down, looseness=0.7] (-0.7,1.4)
to [out=up, in=down] (-0.7,2.4)
to [out=up, in=down, looseness=1.5] (-0.5,2.8);
\strand [thick,-{>[scale=1.7,length=2,width=2]}] (-0.5,0.6)
to [out=up, in=down] (-0.5,2.4)
to [out=up, in=down, looseness=1.5] (-0.3,2.8);
\node[smallnode] (right)   at (0.9,1.6) {};
\strand [thick] (-0.3,0.6) node[below]{$1$}
to [out=up, in=down] (-0.3,2.4)
to [out=up, in=up, looseness=1.5] (0.1,2.4)
to [out=down, in = up] (0.1,1.8)
to [out=down, in = left, looseness = 0.8] (0.9,1.6);
\strand [thick] (0.5,0.6)
to [out=up, in=down] (0.5,2);
\strand [thick,-{>[scale=1.7,length=2,width=2]}] (0.5,2.4)
to [out=up, in=down, looseness=1.5] (0.5,2.8);

\node at (1.6,1.6) {$=$};
\draw[draw=black] (4.3,2) rectangle ++(0.4,0.4);
\node[smallnode] (left)   at (2.2,1.6) {};
\strand [thick,-{>[scale=1.7,length=2,width=2]}] (2.2,1.6)
to [out=right, in=down, looseness=1.5] (2.7,2)
to [out=up, in=down] (2.7,2.4)
to [out=up, in=down, looseness=1.5] (2.9,2.8);
\node at (2.9,0.3) {$n$};
\strand [thick,-{>[scale=1.7,length=2,width=2]}] (2.9,0.6)
to [out=up, in=down] (2.9,2.4)
to [out=up, in=down, looseness=1.5] (3.1,2.8);
\strand [thick,-{>[scale=1.7,length=2,width=2]}] (3.1,0.6)
to [out=up, in=down] (3.1,2.4)
to [out=up, in=down, looseness=1.5] (3.3,2.8);
\strand [thick,-{>[scale=1.7,length=2,width=2]}] (3.3,0.6) node[below]{$j$}
to [out=up,in=down, looseness=0.7] (4.7,1)
to [out=up, in=right, looseness=0.6] (3.6,1.2)
to [out=left, in=down, looseness=1] (3.3,1.4)
to [out=up, in=down] (3.3,2.4)
to [out=up, in=down, looseness=1.5] (3.5,2.8);
\strand [thick,-{>[scale=1.7,length=2,width=2]}] (3.5,0.6)
to [out=up, in=down] (3.5,2.4)
to [out=up, in=down, looseness=1.5] (3.7,2.8);
\node[smallnode] (right)   at (4.9,1.6) {};
\strand [thick] (3.7,0.6) node[below]{$1$}
to [out=up, in=down] (3.7,2.4)
to [out=up, in=up, looseness=1.5] (4.1,2.4)
to [out=down, in = up] (4.1,1.8)
to [out=down, in = left, looseness = 0.8] (4.9,1.6);
\strand [thick] (4.5,0.6)
to [out=up, in=down] (4.5,2);
\strand [thick,-{>[scale=1.7,length=2,width=2]}] (4.5,2.4)
to [out=up, in=down, looseness=1.5] (4.5,2.8);

\node at (5.6,1.6) {$+$};
\draw[draw=black] (8.3,2) rectangle ++(0.4,0.4);
\node[smallnode] (left)   at (6.2,1.6) {};
\strand [thick,-{>[scale=1.7,length=2,width=2]}] (6.2,1.6) node[left]{$z$}
to [out=right, in=down, looseness=1.5] (6.7,2)
to [out=up, in=down] (6.7,2.4)
to [out=up, in=down, looseness=1.5] (6.9,2.8);
\node at (6.9,0.3) {$n$};
\strand [thick,-{>[scale=1.7,length=2,width=2]}] (6.9,0.6)
to [out=up, in=down] (6.9,2.4)
to [out=up, in=down, looseness=1.5] (7.1,2.8);
\strand [thick,-{>[scale=1.7,length=2,width=2]}] (7.1,0.6)
to [out=up, in=down] (7.1,2.4)
to [out=up, in=down, looseness=1.5] (7.3,2.8);
\strand [thick,-{>[scale=1.7,length=2,width=2]}] (7.7,0.6) node[below]{$1$}
to [out=up, in=down] (7.7,1.0)
to [out=up, in=down, looseness=1.8] (7.3,1.4)
to [out=up, in=down] (7.3,2.4)
to [out=up, in=down, looseness=1] (7.5,2.8);
\strand [thick,-{>[scale=1.7,length=2,width=2]}] (7.5,0.6)
to [out=up, in=down] (7.5,2.4)
to [out=up, in=down, looseness=1.5] (7.7,2.8);
\node[smallnode] (right)   at (8.9,1.6) {};
\strand [thick] (7.3,0.6) node[below]{$j$}
to [out=up,in=down, looseness=0.7] (8.7,1)
to [out=up, in=right, looseness=0.7] (7.9,1.2)
to [out=left, in=down] (7.7,1.4)
to [out=up, in=down] (7.7,2.4)
to [out=up, in=up, looseness=1.5] (8.1,2.4)
to [out=down, in = up] (8.1,1.8)
to [out=down, in = left, looseness = 0.8] (8.9,1.6);
\strand [thick] (8.5,0.6)
to [out=up, in=down] (8.5,2);
\strand [thick,-{>[scale=1.7,length=2,width=2]}] (8.5,2.4)
to [out=up, in=down, looseness=1.5] (8.5,2.8);

\flipcrossings{1,2,4,6,7,8,9,13,14,15,16,17,19,20}
\redraw{13}{(3.5,1.1)}
\end{knot}
\draw[thick,blue,dotted] (-0.3,1.2) circle (0.15);
\end{tikzpicture}
\end{equation}

Note that the first term on the right hand side of (\ref{eq:lhs of V_n(1_{n+1})=0}) is equal to $a^{-1} LHS_{j-1}$ via a framed isotopy.\\ 

On the other hand, the right hand side of (\ref{eq:V_n(1_{n+1})=0}) is equal to
\begin{equation}\label{eq:rhs of V_n(1_{n+1})=0}
\begin{tikzpicture}
\tikzstyle{smallnode}=[circle, inner sep=0mm, outer sep=0mm, minimum size=0.7mm, draw=black, fill=black];
\begin{knot}[
end tolerance=.1pt
clip width=5,
clip radius=8pt,
flip crossing/.list={1,2,10,11,12,13,20,21,22,23,24,25},
only when rendering/.style={
}
]
\draw[draw=black] (4.3,2) rectangle ++(0.4,0.4);
\node[smallnode] (left)   at (2.2,1.6) {};
\strand [thick,-{>[scale=1.7,length=2,width=2]}] (2.2,1.6) node[left]{$RHS_j \coloneqq$}
to [out=right, in=left] (2.8,1.6)
to [out=right, in=down, looseness=1.5] (3.3,2)
to [out=up, in=down] (3.3,2.4)
to [out=up, in=down, looseness=1.5] (3.5,2.8);
\node at (2.9,0.3) {$n$};
\strand [thick,-{>[scale=1.7,length=2,width=2]}] (2.9,0.6)
to [out=up, in=down] (2.9,2.4)
to [out=up, in=down, looseness=1.5] (3.1,2.8);
\strand [thick,-{>[scale=1.7,length=2,width=2]}] (3.1,0.6)
to [out=up, in=down] (3.1,2.4)
to [out=up, in=down, looseness=1.5] (3.3,2.8);
\node[smallnode] (right)   at (4.9,1.6) {};
\strand [thick] (3.3,0.6) node[below]{$j$}
to [out=up,in=down] (3.3,1.2)
to [out=up, in=left, looseness=1.5] (3.8,1.6)
to [out=right, in=left] (4.9,1.6);
\strand [thick,-{>[scale=1.7,length=2,width=2]}] (3.5,0.6)
to [out=up, in=down] (3.5,2.4)
to [out=up, in=down, looseness=1.5] (3.7,2.8);
\strand [thick] (3.7,0.6) node[below]{$1$}
to [out=up, in=down] (3.7,2.4)
to [out=up, in=right, looseness=1.5] (3.45,2.6)
to [out=left, in=right] (2.7,2.6)
to [out=left, in=up,looseness=1.5] (2.4,2.3);
\strand [thick,-{>[scale=1.7,length=2,width=2]}] (2.4,2.3)
to [out=down, in=up] (2.4,2)
to [out=down, in=down, looseness=1.5] (2.7,2)
to [out=up, in=down] (2.7,2.4)
to [out=up, in=down, looseness=1.5] (2.9,2.8);
\strand [thick] (4.5,0.6)
to [out=up, in=down] (4.5,2);
\strand [thick,-{>[scale=1.7,length=2,width=2]}] (4.5,2.4)
to [out=up, in=down, looseness=1.5] (4.5,2.8);

\node at (5.6,1.6) {$=$};
\draw[draw=black] (8.3,2) rectangle ++(0.4,0.4);
\node[smallnode] (left)   at (6.2,1.6) {};
\strand [thick,-{>[scale=1.7,length=2,width=2]}] (6.2,1.6) 
to [out=right, in=left] (6.8,1.6)
to [out=right, in=down, looseness=1.5] (7.3,2)
to [out=up, in=down] (7.3,2.4)
to [out=up, in=down, looseness=1.5] (7.5,2.8);
\node at (6.9,0.3) {$n$};
\strand [thick,-{>[scale=1.7,length=2,width=2]}] (6.9,0.6)
to [out=up, in=down] (6.9,2.4)
to [out=up, in=down, looseness=1.5] (7.1,2.8);
\strand [thick,-{>[scale=1.7,length=2,width=2]}] (7.1,0.6)
to [out=up, in=down] (7.1,2.4)
to [out=up, in=down, looseness=1.5] (7.3,2.8);
\node[smallnode] (right)   at (8.9,1.6) {};
\strand [thick] (7.3,0.6) node[below]{$j$}
to [out=up,in=down] (7.3,1.2)
to [out=up, in=left, looseness=1.5] (7.8,1.6)
to [out=right, in=left] (8.9,1.6);
\strand [thick,-{>[scale=1.7,length=2,width=2]}] (7.5,0.6)
to [out=up, in=down] (7.5,2.4)
to [out=up, in=down, looseness=1.5] (7.7,2.8);
\strand [thick] (7.7,0.6) node[below]{$1$}
to [out=up, in=down] (7.7,2.4)
to [out=up, in=right, looseness=1.5] (7.45,2.6)
to [out=left, in=right] (6.7,2.6)
to [out=left, in=up,looseness=1.5] (6.4,2.3);
\strand [thick,-{>[scale=1.7,length=2,width=2]}] (6.4,2.3)
to [out=down, in=up] (6.4,2)
to [out=down, in=down, looseness=1.5] (6.7,2)
to [out=up, in=down] (6.7,2.4)
to [out=up, in=down, looseness=1.5] (6.9,2.8);
\strand [thick] (8.5,0.6)
to [out=up, in=down] (8.5,2);
\strand [thick,-{>[scale=1.7,length=2,width=2]}] (8.5,2.4)
to [out=up, in=down, looseness=1.5] (8.5,2.8);

\node at (9.6,1.6) {$+$};
\draw[draw=black] (12.3,2) rectangle ++(0.4,0.4);
\node[smallnode] (left)   at (10.2,1.6) {};
\strand [thick] (10.2,1.6) node[left]{$z$}
to [out=right, in=left] (10.8,1.6)
to [out=right, in=down, looseness=1.5] (11.3,2)
to [out=up, in=down] (11.3,2.4)
to [out=up, in=right, looseness=1.5] (11.05,2.6)
to [out=left, in=right] (10.7,2.6)
to [out=left, in=up,looseness=1.5] (10.4,2.3);
\strand [thick,-{>[scale=1.7,length=2,width=2]}] (10.4,2.3)
to [out=down, in=up] (10.4,2)
to [out=down, in=down, looseness=1.5] (10.7,2)
to [out=up, in=down] (10.7,2.4)
to [out=up, in=down, looseness=1.5] (10.9,2.8);
\node at (10.9,0.3) {$n$};
\strand [thick,-{>[scale=1.7,length=2,width=2]}] (10.9,0.6)
to [out=up, in=down] (10.9,2.4)
to [out=up, in=down, looseness=1.5] (11.1,2.8);
\strand [thick,-{>[scale=1.7,length=2,width=2]}] (11.1,0.6)
to [out=up, in=down] (11.1,2.4)
to [out=up, in=down, looseness=1.5] (11.3,2.8);
\node[smallnode] (right)   at (12.9,1.6) {};
\strand [thick] (11.3,0.6) node[below]{$j$}
to [out=up,in=down] (11.3,1.2)
to [out=up, in=left, looseness=1.5] (11.8,1.6)
to [out=right, in=left] (12.9,1.6);
\strand [thick,-{>[scale=1.7,length=2,width=2]}] (11.5,0.6)
to [out=up, in=down] (11.5,2.4)
to [out=up, in=down, looseness=1.5] (11.7,2.8);
\strand [thick,-{>[scale=1.7,length=2,width=2]}] (11.7,0.6) node[below]{$1$}
to [out=up, in=down] (11.7,2.4)
to [out=up, in=down, looseness=1.5] (11.5,2.8);
\strand [thick] (12.5,0.6)
to [out=up, in=down] (12.5,2);
\strand [thick,-{>[scale=1.7,length=2,width=2]}] (12.5,2.4)
to [out=up, in=down, looseness=1.5] (12.5,2.8);
\end{knot}
\draw[thick,blue,dotted] (3.4,2.6) circle (0.14);
\end{tikzpicture}
\end{equation}

The first term on the right hand side of (\ref{eq:rhs of V_n(1_{n+1})=0}) is equal to $RHS_{j-1}$ via a framed isotopy, and the second term is equal to $a$ times the second term on the right hand side of (\ref{eq:lhs of V_n(1_{n+1})=0}) via a framed isotopy. This finishes the proof by induction.
\eproof\\

Since $W_n(1_{n+1}) = c^n$, Theorem \ref{thm:relative recursion relation with solution} and Proposition \ref{prop:kernel of V_n} immediately imply the following.

\begin{cor}\label{cor:relative recursion for c^i}
	For every $n \geq 1$ let $A_n \in End_\Lambda(\mathcal{M})$ be so that $[A_i,A_j] = 0$ for all $i,j \in \N$. Then for any $m \in \mathcal{M}$, the relative recursion relation
	\begin{equation}
	[\cdot,e] = r(A, \cdot) \coloneqq r\left(\sum_{i=1}^\infty c^i \otimes A_i, \cdot \right)
	\end{equation}
	has a solution in $\cmm$. \eproof\\
\end{cor}

\begin{ex}\label{ex:recursion and exp/ln}
	Assume that $\mathcal{M} = \left(\mathcal{C}^\Lambda\right)^{\otimes (l-1)}$ for some $l \geq 1$, and that each $A_i, i \in \N,$ is given by multiplication by some element $A_i \in \mathcal{M}$. Since $\mathcal{C}^\Lambda$ is a commutative $\Lambda$-algebra, it follows from Corollary \ref{cor:relative recursion for c^i} that for any $m \in \mathcal{M}$ the relative recursion relation 
	\begin{equation}\label{eq:ex3.7 recursion}
	[\cdot,e] = r\left(\sum_{i=1}^\infty c^i \otimes A_i, \cdot \right)
	\end{equation}
	has a solution in $\cmm$. In the case that $m \in \left(\mathcal{C}_1^\Lambda\right)^{\otimes (l-1)}$, we can even give an explicit formula for the solution $\Phi$: via taking the logarithm/exponential, the equation
	\begin{equation}\label{eq:ex3.7 solved recursion}
		[\Phi,e] = r\left(\sum_{i=1}^\infty c^i \otimes A_i, \Phi \right)
	\end{equation}
	is equivalent to
	\begin{equation}\label{eq:ex3.7 solved recursion after logarithm}
		[\ln \Phi,e] = \sum_{i=1}^\infty c^i \otimes B_i,
	\end{equation}
	where $B_i$ is determined via 
	\begin{equation}
		\sum_{i=1}^\infty c^i \otimes B_i = \ln \left(e \otimes 1 + \sum_{i=1}^\infty c^i \otimes A_i \right).
	\end{equation}

	Indeed, (\ref{eq:ex3.7 solved recursion}) is equivalent to 
	\begin{equation}
		l(\Phi,e) = r\left(e \otimes 1 + \sum_{i=1}^\infty c^i \otimes A_i,\Phi \right) = r(e, \Phi)\cdot \left(e \otimes 1 + \sum_{i=1}^\infty c^i \otimes A_i\right),
	\end{equation}
	which after applying the logarithm becomes
	\begin{equation}
		l(\ln \Phi,e) = r(e, \ln \Phi) + \ln \left(e \otimes 1+ \sum_{i=1}^\infty c^i \otimes A_i\right) = r(e, \ln \Phi) + \sum_{i=1}^\infty c^i \otimes B_i,
	\end{equation}
	which is equivalent to (\ref{eq:ex3.7 solved recursion after logarithm}). The degree zero part of $\ln \Phi$ is given by\footnote{This is where we use that $m \in \left(\mathcal{C}_1^\Lambda\right)^{\otimes (l-1)}$ so that $\ln(m)$ is defined.}
	\begin{equation}
		\ln(1 \otimes m) = 1 \otimes \ln(m).
	\end{equation} 
	Furthermore, Morton \cite{mor02} showed that the elements $P_i \in \mathcal{C}^+$ defined as 
	\begin{equation}\label{eq:def of P_i}
		\begin{tikzpicture}
			\tikzstyle{smallnode}=[circle, inner sep=0mm, outer sep=0mm, minimum size=0.5mm, draw=black, fill=black];
			\begin{knot}[
				ignore endpoint intersections=false,
				clip width=5,
				clip radius=8pt,
				only when rendering/.style={
				}
				]
				\node at (-1.2,0.5) {$P_i=\frac{s-s^{-1}}{s^i-s^{-i}}\Bigg($};
				\node at (9.3,0.5) {$\Bigg)$};
				
				\strand [thick,-{>[scale=2,length=2,width=3]}] (0,0)
				to (1.4,1);
				\strand [thick,-{>[scale=2,length=2,width=3]}] (0.2,0)
				to (0,1);
				\strand [thick,-{>[scale=2,length=2,width=3]}] (0.4,0)
				to (0.2,1);
				\strand [thick,-{>[scale=2,length=2,width=3]}] (0.601,0)
				to (0.4,1);
				\strand [thick,-{>[scale=2,length=2,width=3]}] (1.4,0)
				to (1.2,1);
				\node at (1,0.3) {$...$};
				\node at (0,-0.3) {$i$};
				\node at (1.4,-0.3) {$1$};
				
				\node at (1.8,0.5) {$+$};
				\strand [thick,-{>[scale=2,length=2,width=3]}] (2.2,0)
				to (3.6,1);
				\strand [thick,-{>[scale=2,length=2,width=3]}] (2.4,0)
				to (2.2,1);
				\strand [thick,-{>[scale=2,length=2,width=3]}] (2.601,0)
				to (2.4,1);
				\strand [thick,-{>[scale=2,length=2,width=3]}] (2.8,0)
				to (2.6,1);
				\strand [thick,-{>[scale=2,length=2,width=3]}] (3.6,0)
				to (3.4,1);
				\node at (3.2,0.3) {$...$};
				
				\node at (4,0.5) {$+$};
				\strand [thick,-{>[scale=2,length=2,width=3]}] (4.4,0)
				to (5.8,1);
				\strand [thick,-{>[scale=2,length=2,width=3]}] (4.6,0)
				to (4.4,1);
				\strand [thick,-{>[scale=2,length=2,width=3]}] (4.801,0)
				to (4.6,1);
				\strand [thick,-{>[scale=2,length=2,width=3]}] (5.0,0)
				to (4.8,1);
				\strand [thick,-{>[scale=2,length=2,width=3]}] (5.8,0)
				to (5.6,1);
				\node at (5.4,0.3) {$...$};
				
				\node at (6.7,0.5) {+ ... +};
				\strand [thick,-{>[scale=2,length=2,width=3]}] (7.6,0)
				to (9,1);
				\strand [thick,-{>[scale=2,length=2,width=3]}] (7.8,0)
				to (7.6,1);
				\strand [thick,-{>[scale=2,length=2,width=3]}] (8,0)
				to (7.8,1);
				\strand [thick,-{>[scale=2,length=2,width=3]}] (8.2,0)
				to (8,1);
				\strand [thick,-{>[scale=2,length=2,width=3]}] (9,0)
				to (8.8,1);
				\node at (8.6,0.3) {$...$};		
				
				\flipcrossings{5,9,10,13,14,15,16}			
				\redraw{7}{(2.35,0.1)}
				\redraw{13}{(4.75,0.2)}
				\redraw{20}{(8.05,0.1)}
			\end{knot}
		\end{tikzpicture},
	\end{equation}
	which correspond to the power sums $\sum_k x_k^i$ under the isomorphism of $\mathcal{C}^+$ and the ring of symmetric functions by \cite[Theorem 4.9.6]{ais96}\footnote{A small warning: what Aiston calls $P_m$ corresponds to $\frac{s^m-s^{-m}}{s-s^{-1}}P_m$ in our language.}, satisfy the commutator relation 
	\begin{equation}\label{eq:relative commutator of P_i}
		[P_i,e] = (s^{i}-s^{-i}) c^i.
	\end{equation}
	It follows that $\ln \Phi$ is given by
	\begin{equation}\label{eq:ex3.7 logaritmic solution}
		\ln \Phi = 1 \otimes \ln(m) + \sum_{i=1}^\infty \frac{1}{s^i-s^{-i}} P_i \otimes B_i.
	\end{equation}
	By exponentiating (\ref{eq:ex3.7 logaritmic solution}), we obtain a formula for $\Phi$.
	
	Conversely, if $\Phi \in \cmm$ satisfies 
	\begin{equation}
		[\Phi,e] = \sum_{i=1}^\infty c^i \otimes B_i,
	\end{equation}
	where $m \in \left(\mathcal{C}_0^\Lambda\right)^{\otimes (l-1)}$, then $\exp \Phi \in \mathcal{C}^{\mathcal{M}}_{\exp m}$ satisfies
	\begin{equation}
		[\exp \Phi,e] = r\left(\sum_{i=1}^\infty c^i \otimes A_i, \exp \Phi \right),
	\end{equation}
	where $A_i$ is determined via 
	\begin{equation}
		\sum_{i=1}^\infty c^i \otimes A_i =  \exp \left(\sum_{i=1}^\infty c^i \otimes B_i \right) - e \otimes 1.
	\end{equation}\\
\end{ex}

\begin{ex}
Consider the partition function 
\begin{equation}
H(t) = 1 + \sum\limits_{i = 1}^\infty W_{(i)} t^i,
\end{equation}
which encodes the basis elements $W_\lambda$ for all partitions $\lambda = (l(\lambda))$ whose Young diagram consists of a single row. Its logarithm is given by
\begin{equation}
\ln H(t) = \sum\limits_{i=1}^\infty \frac{1}{i} P_i t^i.
\end{equation}

As was already observed by Pokorny \cite[Theorem 4.1.4]{pok21}, this implies that $H(t)$ satisfies the recursion relation 
\begin{equation}
[H(t),e] = r(A(t),H(t)),
\end{equation}
where $A(t)$ is given by 
\begin{equation}
A(t) = \exp\left(\sum_{i=1}^\infty \frac{s^i- s^{-i}}{i} c^i t^i\right) - e = \frac{1-s^{-1}ct}{1-sct} - e = \left(1-s^{-2} \right) \sum\limits_{i=1}^\infty s^i c^i t^i.  
\end{equation}\\
\end{ex}

\begin{rmk}
	By performing the mirror of the above arguments, one obtains analogous results for relative recursion relations with $l$ instead of $r$.\\
\end{rmk}

\section{The BPS partition functions $\Psi^{(g,l)}$.}\label{sec:bps partition functions}

In this section, we introduce the BPS partition functions and prove Theorem \ref{thm:properties of bps partition functions}.\\

Using the power sum elements $P_i$ defined in (\ref{eq:def of P_i}), we define the partition functions
\begin{equation}
\Psi^{(g,l)} \coloneqq \exp\left( \sum\limits_{i=1}^\infty \frac{1}{i} \frac{1}{(s^i - s^{-i})^{2-2g-l}} (P_i)^{\otimes l}  \right)
\end{equation}
with values in $\mathcal{C}^{\otimes l}$ for any $l \in \N_0$ and $g \in \Z$. Their inverses in the skein are given by 
\begin{equation}
\Psi^{-(g,l)} \coloneqq  \exp\left( -\sum\limits_{i=1}^\infty \frac{1}{i} \frac{1}{(s^i - s^{-i})^{2-2g-l}} (P_i)^{\otimes l}  \right).
\end{equation}

Here, we set $\mathcal{C}^{\otimes l} = \Lambda$ and $(P_i)^{\otimes l} = 1$ if $l = 0$. Note that these partition functions are invariant under interchanging the different copies of the solid torus.\\

We will denote the logarithm of these partition functions by 
\begin{equation}
\psi^{\pm(g,l)} \coloneqq \ln \Psi^{\pm(g,l)} = \pm \sum\limits_{i=1}^\infty \frac{1}{i} \frac{1}{(s^i - s^{-i})^{2-2g-l}} (P_i)^{\otimes l}.
\end{equation}\\

Clearly, we have that $\Psi^{(g,l)} \Psi^{-(g,l)} = 1$.

Observe that $\Psi^{(g,l)}$ is invariant under the mirror map if $l$ is even, and that the mirror of $\Psi^{(g,l)}$ is $\Psi^{-(g,l)}$ if $l$ is odd since $P_i$ is invariant under the mirror map by \cite[p.\,315]{mm08}.

\subsection{Crossing a Lagrangian}

Let $\mathcal{U}: \mathcal{C} \to \mathcal{S}(S^3) \approx \Lambda$ be the $\Lambda$-algebra homomorphism given by decorating the unknot in $S^3$ with the standard (zero) framing with an element of $\mathcal{C}$. Then the following holds.
\begin{lem}\label{lem:Psi^(g,l-1) crossing a lagrangian} For any $l \geq 1$, we have the identity
\begin{equation}\label{eq:psi^(g,l-1) crossing a lagrangian}
\Psi^{\pm(g,l-1)}(a t) \Psi^{\mp(g,l-1)}(a^{-1} t) = \left( 1^{\otimes l-1} \otimes \mathcal{U} \right)(\Psi^{\pm(g,l)}).
\end{equation}
\end{lem}

\bproof Since $\mathcal{U}$ is an algebra homomorphism, the logarithm commutes with $1^{\otimes l-1} \otimes \mathcal{U}$. Thus, via taking the logarithm on both sides, (\ref{eq:psi^(g,l-1) crossing a lagrangian}) is equivalent to 
\begin{equation}
\sum\limits_{i=1}^\infty \pm \frac{1}{i} \frac{a^i - a^{-i}}{(s^i - s^{-i})^{2-2g-(l-1)}} (P_i)^{\otimes l-1} =  \sum\limits_{i=1}^\infty \pm \frac{1}{i} \frac{1}{(s^i - s^{-i})^{2-2g-l}} (P_i)^{\otimes l-1} \mathcal{U}(P_i),
\end{equation}
which is satisfied since $\mathcal{U}(P_i) = \frac{a^i-a^{-i}}{s^i-s^{-i}}$ which follows directly from the explicit representation of $P_i$ given in (\ref{eq:def of P_i}) above. \eproof\\

\subsection{The partition function $\Psi^{(1,1)}$ of a torus with one disk removed.} Let us consider the case $(g=1,l=1)$. It was shown by Morton and Manch\'{o}n \cite[Theorem 1, Theorem 2]{mm08}, that
\begin{equation}\label{eq:geometric expression for Psi^(1,1)}
\Psi^{(1,1)} = 1 + z \sum_{i=1}^\infty C_i
\end{equation}
and 
\begin{equation}\label{eq:geometric expression for Psi^-(1,1)}
\Psi^{-(1,1)} = 1 - z \sum_{i=1}^\infty \overline{C}_i,
\end{equation}
where $\overline{C}_i$ denotes the mirror of $C_i$.\\

\bproof
We will give a proof which is different from the proof in \cite{mm08}, using recursion relations. By (\ref{eq:relative commutator of P_i}), we have the commutation relation 
\begin{equation}
[\psi^{(1,1)},e] = \sum_{i=1}^\infty \frac{(s^i - s^{-i})^2}{i} c^i.
\end{equation}

By exponentiating (see Example \ref{ex:recursion and exp/ln} above), we obtain the recursion relation
\begin{equation}
[\Psi^{(1,1)},e] = r\left( \frac{(1-c)^2}{(1-s^2c)(1-s^{-2}c)} - e, \Psi^{(1,1)} \right),
\end{equation}
where
\begin{equation}
\begin{split}
\frac{(1-c)^2}{(1-s^2c)(1-s^{-2}c)} - e&= (1-c)^2 \left( \sum_{i=0}^\infty s^{2i} c^i \right) \left( \sum_{i=0}^\infty s^{-2i} c^i \right) - e\\
&= \sum_{i=1}^\infty (s^{2i} - 2s^{2i-2} + 2s^{2i-4} - 2s^{2i-6} + ...- 2s^{2-2i} + s^{-2i}) c^i\\
&= \sum_{i=1}^\infty (s-s^{-1}) \frac{s^{2i} - s^{-2i}}{s + s^{-1}} c^i.
\end{split}
\end{equation}

We show that $1 + z \sum_{i=1}^\infty C_i$ satisfies the same recursion relation. First observe that 
\begin{equation}
[C_i,e] = \left[a^{-1} \frac{P - \bigcirc}{s^i - s^{-i}} (P_i),e \right] = a^{-1} (P - \bigcirc) c^i
\end{equation}
since $P - \bigcirc$ and $[\cdot, e]$ ``commute", where we used that 
\begin{equation}
	C_i = a^{-1} \frac{P - \bigcirc}{s^i - s^{-i}} (P_i),
\end{equation}
which follows from (\ref{eq:p-unknot on e_lambda}) above and the expansions of $C_i$ and $P_i$ in terms of the $W_\lambda$ given in \cite[Theorem 15]{mm08} and its preceding paragraphs.

Now we show that
\begin{equation}\label{eq:recursion for C(t)}
a^{-1} (P - \bigcirc) c^i = \sum_{j=1}^i z \frac{s^{2j} - s^{-2j}}{s + s^{-1}} r(c^j, C_{i-j})
\end{equation}
by induction on $i$, where we have set $C_0 = \frac{1}{z}$. 

For $i = 1$, (\ref{eq:recursion for C(t)}) follows immediately from the skein relations. Assume that (\ref{eq:recursion for C(t)}) is true for all $i < n $ for some $n \geq 2$. By applying the skein relations, we see that
\begin{equation}
	\begin{tikzpicture}
		\tikzstyle{smallnode}=[circle, inner sep=0mm, outer sep=0mm, minimum size=0.7mm, draw=black, fill=black];
		\begin{knot}[
			ignore endpoint intersections=false,
			clip width=5,
			clip radius=8pt,
			only when rendering/.style={
			},
			]
			\node[smallnode] at (-0.4,0.5) {};
			\strand [thick,looseness=1,-{>[scale=2,length=2,width=3]}] (-0.4,0.5)
			to [out=0, in=-90] (0.4,2); 
			\strand [thick,-{>[scale=2,length=2,width=3]}] (0.4,0)
			to [out=90, in=-90] (0.8,2); 
			\strand [thick,-{>[scale=2,length=2,width=3]}] (0.8,0)
			to [out=90, in=-90] (1.2,2); 
			\strand [thick,-{>[scale=2,length=2,width=3]}] (1.2,0)
			to [out=90, in=-90] (1.6,2); 
			\strand [thick,-{>[scale=2,length=2,width=3]}] (1.6,0)
			to [out=90, in=-90] (2,2); 
			\strand [thick,-{>[scale=2,length=2,width=3]}] (2,0)
			to [out=90, in=-90] (2.4,2);
			\node[smallnode] at (3.2,0.5) {};
			\strand [thick,looseness=1.4] (2.4,0)
			to [out=90, in=180] (3.2,0.5); 
			\strand [thick,looseness=0.5,-{>[scale=2,length=1.5,width=5]}] (0,1.4)
			to [out=-90,in=180] (1.4,1.1)
			to [out=0,in=-90] (2.8,1.4)
			to [out=90,in=0] (1.4,1.7)
			to [out=180,in=90] (0,1.4);
			\node at (0.4,2.3) {$n$};
			\node at (2.4,2.3) {$1$};
			
			\node at (3.8,1) {$=$};
			
			\node[smallnode] at (4.4,0.5) {};
			\strand [thick,looseness=1,-{>[scale=2,length=2,width=3]}] (4.4,0.5)
			to [out=0, in=-90] (5.2,2); 
			\strand [thick,-{>[scale=2,length=2,width=3]}] (5.2,0)
			to [out=90, in=-90] (5.6,2); 
			\strand [thick,-{>[scale=2,length=2,width=3]}] (5.6,0)
			to [out=90, in=-90] (6,2); 
			\strand [thick,-{>[scale=2,length=2,width=3]}] (6,0)
			to [out=90, in=-90] (6.4,2); 
			\strand [thick,-{>[scale=2,length=2,width=3]}] (6.4,0)
			to [out=90, in=-90] (6.8,2); 
			\strand [thick,-{>[scale=2,length=2,width=3]}] (6.8,0)
			to [out=90, in=-90] (7.2,2);
			\node[smallnode] at (8,0.5) {};
			\strand [thick,looseness=1.4] (7.2,0)
			to [out=90, in=180] (8,0.5); 
			\strand [thick,looseness=0.5,-{>[scale=2,length=1.5,width=5]}] (4.8,1.4)
			to [out=-90,in=180] (5.9,1.1)
			to [out=0,in=-90] (7.0,1.4)
			to [out=90,in=0] (5.9,1.7)
			to [out=180,in=90] (4.8,1.4);
			\node at (5.2,2.3) {$n$};
			\node at (7.2,2.3) {$1$};
			
			\node at (8.6,1) {$+z$}; 
			
			\node[smallnode] at (9.2,0.5) {};
			\strand [thick,looseness=1,-{>[scale=2,length=2,width=3]}] (9.2,0.5)
			to [out=0, in=-90] (10,2); 
			\strand [thick,-{>[scale=2,length=2,width=3]}] (10,0)
			to [out=90, in=-90] (10.4,2); 
			\strand [thick,-{>[scale=2,length=2,width=3]}] (10.4,0)
			to [out=90, in=-90] (10.8,2); 
			\strand [thick,-{>[scale=2,length=2,width=3]}] (10.8,0)
			to [out=90, in=-90] (11.2,2); 
			\strand [thick,-{>[scale=2,length=2,width=3]}] (11.2,0)
			to [out=90, in=-90] (11.6,2); 
			\strand [thick,looseness=1] (11.6,0)
			to [out=90,in=250] (11.85,1.15)
			to [out=70,in=-90] (12.4,1.4);
			\strand [thick,looseness=0.5] (12.4,1.4)
			to [out=90,in=0] (11,1.7)
			to [out=180,in=90] (9.6,1.4)
			to [out=-90,in=180] (10.7,1.1)
			to [out=0,in=200] (11.75,1.2);
			\strand [thick,looseness=0.5,-{>[scale=2,length=2,width=3]}] (11.75,1.2)
			to [out=20, in=-90] (12,2);
			\node[smallnode] at (12.8,0.5) {};
			\strand [thick,looseness=1.4] (12,0)
			to [out=90, in=180] (12.8,0.5); 
			\node at (10,2.3) {$n$};
			\node at (12,2.3) {$1$};
			
			\flipcrossings{1,3,5,7,9,11,13,15,17,19,21,24,26,28,30,32,34}
		\end{knot}
	\end{tikzpicture}
\end{equation}
and 
\begin{equation}
	\begin{tikzpicture}
		\tikzstyle{smallnode}=[circle, inner sep=0mm, outer sep=0mm, minimum size=0.7mm, draw=black, fill=black];
		\begin{knot}[
			ignore endpoint intersections=false,
			clip width=5,
			clip radius=8pt,
			only when rendering/.style={
			},
			]
			\node[smallnode] at (-0.4,0.5) {};
			\strand [thick,looseness=1,-{>[scale=2,length=2,width=3]}] (-0.4,0.5)
			to [out=0, in=-90] (0.4,2); 
			\strand [thick,-{>[scale=2,length=2,width=3]}] (0.4,0)
			to [out=90, in=-90] (0.8,2); 
			\strand [thick,-{>[scale=2,length=2,width=3]}] (0.8,0)
			to [out=90, in=-90] (1.2,2); 
			\strand [thick,-{>[scale=2,length=2,width=3]}] (1.2,0)
			to [out=90, in=-90] (1.6,2); 
			\strand [thick,-{>[scale=2,length=2,width=3]}] (1.6,0)
			to [out=90, in=-90] (2,2); 
			\strand [thick,looseness=1] (2,0)
			to [out=90,in=250] (2.25,1.15)
			to [out=70,in=-90] (2.8,1.4);
			\strand [thick,looseness=0.5] (2.8,1.4)
			to [out=90,in=0] (1.4,1.7)
			to [out=180,in=90] (0,1.4)
			to [out=-90,in=180] (1.1,1.1)
			to [out=0,in=200] (2.15,1.2);
			\strand [thick,looseness=0.5,-{>[scale=2,length=2,width=3]}] (2.15,1.2)
			to [out=20, in=-90] (2.4,2);
			\node[smallnode] at (3.2,0.5) {};
			\strand [thick,looseness=1.4] (2.4,0)
			to [out=90, in=180] (3.2,0.5);
			\node at (0.4,2.3) {$n$};
			\node at (2.4,2.3) {$1$};
			
			\node at (3.8,1) {$=$};
			
			\node[smallnode] at (4.4,0.5) {};
			\strand [thick,looseness=1,-{>[scale=2,length=2,width=3]}] (4.4,0.5)
			to [out=0, in=-90] (5.2,2); 
			\strand [thick,-{>[scale=2,length=2,width=3]}] (5.2,0)
			to [out=90, in=-90] (5.6,2); 
			\strand [thick,-{>[scale=2,length=2,width=3]}] (5.6,0)
			to [out=90, in=-90] (6,2); 
			\strand [thick,-{>[scale=2,length=2,width=3]}] (6,0)
			to [out=90, in=-90] (6.4,2); 
			\strand [thick,-{>[scale=2,length=2,width=3]}] (6.4,0)
			to [out=90, in=-90] (6.8,2); 
			\strand [thick,looseness=1] (6.8,0)
			to [out=90,in=250] (7.05,1.15)
			to [out=70,in=-90] (7.6,1.4);
			\strand [thick,looseness=0.5] (7.6,1.4)
			to [out=90,in=0] (6.2,1.7)
			to [out=180,in=90] (4.8,1.4)
			to [out=-90,in=180] (5.9,1.1)
			to [out=0,in=200] (6.96,1.2);
			\strand [thick,looseness=0.5,-{>[scale=2,length=2,width=3]}] (6.96,1.2)
			to [out=20, in=-90] (7.2,2);
			\node[smallnode] at (8,0.5) {};
			\strand [thick,looseness=1.4] (7.2,0)
			to [out=90, in=180] (8,0.5); 
			\node at (5.2,2.3) {$n$};
			\node at (7.2,2.3) {$1$};
			
			\node at (8.6,1) {$+z \sum\limits_{j=2}^n$}; 
			
			\node[smallnode] at (9.2,0.5) {};
			\strand [thick,looseness=1,-{>[scale=2,length=2,width=3]}] (9.2,0.5)
			to [out=0, in=-90] (10,2); 
			\strand [thick,-{>[scale=2,length=2,width=3]}] (10,0)
			to [out=90, in=-90] (10.4,2); 
			\strand [thick] (10.4,0)
			to [out=90,in=250] (10.6,1);
			\strand [thick,looseness=0.5] (10.6,1)
			to [out=70,in=200] (11.75,1.2);
			\strand [thick,-{>[scale=2,length=2,width=3]},looseness=0.5] (11.75,1.2)
			to [out=20, in=-90] (12,2);
			\strand [thick,-{>[scale=2,length=2,width=3]}] (10.8,0)
			to [out=90, in=-90] (11.2,2); 
			\strand [thick,-{>[scale=2,length=2,width=3]}] (11.2,0)
			to [out=90, in=-90] (11.6,2); 
			\strand [thick,looseness=1] (11.6,0)
			to [out=90,in=250] (11.85,1.15)
			to [out=70,in=-90] (12.4,1.4);
			\strand [thick,looseness=0.5] (12.4,1.4)
			to [out=90,in=0] (11,1.7)
			to [out=180,in=90] (9.6,1.4)
			to [out=-90,in=190] (10.55,1.1);
			\strand [thick,looseness=0.5,-{>[scale=2,length=2,width=3]}] (10.55,1.1)
			to [out=10, in=-90] (10.8,2); 
			\node[smallnode] at (12.8,0.5) {};
			\strand [thick,looseness=1.4] (12,0)
			to [out=90, in=180] (12.8,0.5); 
			\node at (10,2.3) {$n$};
			\node at (12,2.3) {$1$};
			\node at (10.8,2.3) {$j$};
			
			\redraw{14}{(6.56,1.15)}
			\redraw{20}{(10.15,1.1)}
			\flipcrossings{2,4,6,8,10,12,25,39}
		\end{knot}
	\end{tikzpicture},
\end{equation}
which implies that
\begin{equation}
a^{-1} (P - \bigcirc) c^n = c \cdot \left( a^{-1} (P - \bigcirc)\left(c^{n-1}\right) \right) + z c^n + \sum_{j=2}^n z^2 l(C_{n-j+1},c^{j-1}). 
\end{equation}

By the induction hypothesis, we have that
\begin{equation}
a^{-1} (P - \bigcirc) \left(c^{n-1}\right) = \sum_{j=1}^{n-1} z \frac{s^{2j} - s^{-2j}}{s + s^{-1}} r(c^j, C_{n-1-j})
\end{equation}
and
\begin{equation}
\begin{split}
l(C_{n-j+1},c^{j-1}) &= r(c^{j-1},C_{n-j+1}) + c^{j-1} \cdot [C_{n-j+1},e]\\
&= r(c^{j-1},C_{n-j+1}) +  \sum_{k=1}^{n-j+1} z \frac{s^{2k} - s^{-2k}}{s + s^{-1}} r(c^{k+j-1}, C_{n-j+1-k}).
\end{split}
\end{equation}

Thus, $a^{-1} (P - \bigcirc) c^n$ can be written as 
\begin{eqnarray}
a^{-1} (P - \bigcirc) c^n = \sum_{l=1}^n a_l r(c^l, C_{n-l}),
\end{eqnarray}
where the coefficient $a_l \in \Lambda$ is given by 
\begin{equation}
\begin{gathered}
z \frac{s^{2(l-1)} - s^{-2(l-1)}}{s + s^{-1}} + z^2 + \sum_{j=2}^l z^3 \frac{s^{2(l+1-j)} - s^{-2(l+1-j)}}{s + s^{-1}}\\
= \frac{z}{s + s^{-1}} \left( s^{2(l-1)} - s^{-2(l-1)} + s^2 - s^{-2} + z \left(s^{2l-1} - s - s^{-1} + s^{-2l+1} \right) \right) \\
=  \frac{z}{s + s^{-1}} \left( s^{2l} - s^{-2l} \right).
\end{gathered}
\end{equation}

This finishes the induction step and the proof of (\ref{eq:geometric expression for Psi^(1,1)}). The formula (\ref{eq:geometric expression for Psi^-(1,1)}) follows from (\ref{eq:geometric expression for Psi^(1,1)}) by applying the mirror map.
\eproof\\

\subsection{The disk partition function $\Psi^{(0,1)}$}

Let us consider the case $(g=0,l=1)$. By (\ref{eq:relative commutator of P_i}), we have the commutation relation 
\begin{equation}
[\psi^{(0,1)},e] = \sum_{i=1}^\infty \frac{1}{i} c^i.
\end{equation}
As $\exp\left(\sum_{i=1}^\infty \frac{1}{i} c^i \right) = \frac{1}{1-c} = \sum_{i=0}^\infty c^i$, it follows that $\Psi^{(0,1)}$ satisfies the relative recursion relation
\begin{equation}
[\Psi^{(0,1)},e] = r\left(\sum_{i=1}^\infty c^i, \Psi^{(0,1)}\right).
\end{equation}

Note that this implies that 
\begin{equation}\label{eq:variant relative recursion for Psi^(0,1)}
l(\Psi^{(0,1)},c) = r(c, \Psi^{(0,1)}) + c \cdot [\Psi^{(0,1)},e] = r\left(\sum_{i=1}^\infty c^i, \Psi^{(0,1)}\right) = [\Psi^{(0,1)},e].
\end{equation}

By applying the mirror map, we also obtain
\begin{equation}\label{eq:relative recursion for Psi^-(0,1)}
[\Psi^{-(0,1)},e] = -r(c,\Psi^{-(0,1)}).
\end{equation}

By applying the wiring $\cp$ defined in (\ref{eq:definition of cp}), (\ref{eq:relative recursion for Psi^-(0,1)}) becomes $(P_{1,0} - \bigcirc + a P_{0,1})\Psi^{-(0,1)} = 0$ in the notation of Ekholm and Shende \cite{es20}, where it was shown that
\begin{equation}
\Psi^{-(0,1)} = 1 + \sum_{n\geq1} \sum_{\lambda \vdash n} (-1)^{n} W_\lambda \prod_{\square \in \lambda} \frac{s^{-c(\square)}}{s^{h(\square)} - s^{-h(\square)}},
\end{equation}
where $c(\square)$ and $h(\square)$ denote the content and hook length of a box in the young diagram of an (unordered) partition $\lambda$, respectively. By applying the mirror map again, we obtain
\begin{equation}
\Psi^{(0,1)} = 1 + \sum_{n\geq1} \sum_{\lambda \vdash n} W_\lambda \prod_{\square \in \lambda} \frac{s^{c(\square)}}{s^{h(\square)} - s^{-h(\square)}}.
\end{equation}\\

\subsection{The annulus partition function $\Psi^{(0,2)}$}

Let us consider the case $(g=0,l=2)$. By (\ref{eq:relative commutator of P_i}), we have the commutation relation 
\begin{equation}
\left( [\cdot,e] \otimes 1 \right) (\psi^{(0,2)}) = \sum_{i=1}^\infty \frac{1}{i} (s^i - s^{-i}) c^i \otimes P_i.
\end{equation}
It follows that $\Psi^{(0,2)}$ satisfies the relative recursion relation
\begin{equation}\label{eq:relative recursion for Psi^(0,2)}
\left([\cdot,e] \otimes 1\right)(\Psi^{(0,2)}) = r\left(\exp\left(\sum_{i=1}^\infty \frac{1}{i} (s^i - s^{-i}) c^i \otimes P_i \right) - e, \Psi^{(0,2)} \right).
\end{equation}
In order to obtain a nice expression for $\exp\left(\sum_{i=1}^\infty \frac{1}{i} (s^i - s^{-i}) c^i \otimes P_i \right)$, we use that $\Psi^{(1,1)} = \exp\left(\sum_{i=1}^\infty \frac{1}{i} (s^i - s^{-i}) P_i \right) = 1 + z \sum_{i=1}^\infty C_i$. It follows that
\begin{equation}
\exp\left(\sum_{i=1}^\infty \frac{1}{i} (s^i - s^{-i}) c^i \otimes P_i \right) - e = z \sum_{i=1}^\infty c^i \otimes C_i.
\end{equation}

By applying $\cp \otimes 1$, (\ref{eq:relative recursion for Psi^(0,2)}) implies that
\begin{equation}
((P - \bigcirc) \otimes 1)(\Psi^{(0,2)}) = \left(az \sum_{i=1}^\infty C_i \otimes C_i \right) \Psi^{(0,2)}.
\end{equation}

By symmetry, we also have that
\begin{equation}
(1 \otimes (P - \bigcirc))(\Psi^{(0,2)}) = \left(az \sum_{i=1}^\infty C_i \otimes C_i \right) \Psi^{(0,2)} = ((P - \bigcirc) \otimes 1)(\Psi^{(0,2)}).
\end{equation}

Since $P - \bigcirc$ has one-dimensional eigenspaces spanned by the $W_\lambda$'s, it follows that $\Psi^{(0,2)}$ is of the form
\begin{equation}\label{eq:Psi^(0,2) in terms of W_lambda}
\Psi^{(0,2)} = 1 + \sum_{n\geq1} \sum_{\lambda \vdash n} a_\lambda W_\lambda \otimes W_\lambda
\end{equation}
for some coefficients $a_\lambda \in \Lambda$.

We claim that $a_\lambda = 1$ for all partitions $\lambda$. First note that $a_\lambda$ does not depend on $a$. Indeed, the expression of $\psi^{(0,2)}$ in terms of the $P_i$ does not depend on $a$ and any polynomial in the $P_i$'s can be written as a linear combination of $W_\lambda$'s with coefficients in $\Q$. To see that $a_\lambda = 1$, we apply Lemma \ref{lem:Psi^(g,l-1) crossing a lagrangian} and use the fact (see  \cite[Corollary 6.7]{ais97}) that the coefficient of the highest power of $a$ in $\mathcal{U}(W_\lambda)$ is given by
\begin{equation}
\prod_{\square \in \lambda} \frac{s^{c(\square)}}{s^{h(\square)} - s^{-h(\square)}}.
\end{equation}
Thus, comparing the coefficient of the highest powers of $a$ in $(1 \otimes \mathcal{U})(\Psi^{(0,2)})$ and $\Psi^{(0,1)}(at) \Psi^{-(0,1)}(a^{-1}t)$ in each degree gives
\begin{equation}
\sum_\lambda W_\lambda \prod_{\square \in \lambda} \frac{s^{c(\square)}}{s^{h(\square)} - s^{-h(\square)}} = \sum_\lambda a_\lambda W_\lambda \prod_{\square \in \lambda} \frac{s^{c(\square)}}{s^{h(\square)} - s^{-h(\square)}},
\end{equation}
which implies that $a_\lambda = 1$ as claimed.\\

By analogous arguments, it follows that 
\begin{equation}
\left([\cdot,e] \otimes 1\right)(\Psi^{-(0,2)}) = r\left(- z \sum_{i=1}^\infty c^i \otimes \overline{C_i}, \Psi^{(0,2)} \right),
\end{equation}
which implies that 
\begin{equation}
a^{-1} ((P - \bigcirc) \otimes 1)(\Psi^{-(0,2)}) = \left(-z \sum_{i=1}^\infty C_i \otimes \overline{C_i} \right) \Psi^{-(0,2)} = -\overline{a^{-1} (1 \otimes (P - \bigcirc))(\Psi^{-(0,2)})}
\end{equation}
by symmetry. Since 
\begin{equation}
a^{-1} ((P - \bigcirc) \otimes 1)(W_\lambda \otimes W_\mu) = z \sum_{\square \in \lambda} s^{2 c(\square)} W_\lambda \otimes W_\mu
\end{equation}
and 
\begin{equation}
-\overline{a^{-1} (1 \otimes (P - \bigcirc))(W_\lambda \otimes W_\mu)} = z \sum_{\square \in \mu} s^{-2 c(\square)} W_\lambda \otimes W_\mu,
\end{equation}
where we used that $W_\lambda$ is invariant under the mirror map, and because a partition $\lambda$ is uniquely determined by $\sum_{\square \in \lambda} s^{2 c(\square)} = \sum_{\square \in \lambda'} s^{-2 c(\square)}$, where $\lambda'$ denotes the partition conjugate to $\lambda$, it follows that
\begin{equation}
\Psi^{-(0,2)} = 1 + \sum_{\{\lambda,\lambda'\}} a_{\{\lambda,\lambda'\}} (W_\lambda \otimes W_{\lambda'} + W_{\lambda'} \otimes W_\lambda)
\end{equation}
for some coefficients $a_{\{\lambda,\lambda'\}} \in \Lambda$, which as above do not depend on $a$. By comparing the coefficient of the highest powers of $a$ in $(1 \otimes \mathcal{U})(\Psi^{-(0,2)})$ and $\Psi^{-(0,1)}(at) \Psi^{(0,1)}(a^{-1}t)$ in each degree, we obtain that
\begin{equation}
\sum_\lambda (-1)^{l(\lambda)} W_\lambda \prod_{\square \in \lambda} \frac{s^{-c(\square)}}{s^{h(\square)} - s^{-h(\square)}} = \sum_\lambda a_{\{\lambda,\lambda'\}} 2^{\delta_{\lambda,\lambda'}} W_\lambda \prod_{\square \in \lambda'} \frac{s^{c(\square)}}{s^{h(\square)} - s^{-h(\square)}}.
\end{equation}
It follows that $a_{\{\lambda,\lambda'\}} = \frac{(-1)^{l(\lambda)}}{2}$ if $\lambda = \lambda'$ and that $a_{\{\lambda,\lambda'\}} = (-1)^{l(\lambda)}$ if $\lambda \neq \lambda'$. Thus, we obtain the formula
\begin{equation}
\Psi^{-(0,2)} = 1 + \sum_{n \geq 1} (-1)^n \sum_{\lambda \vdash n} W_\lambda \otimes W_{\lambda'}.
\end{equation}\\

\subsection{Gluing formulas for $\Psi^{(g,l)}$} Denote by $W_\lambda^* \in \left( \mathcal{C^+}\right)^*$ the basis elements dual to $W_\lambda$. Similarly, we define $P_\lambda \coloneqq \prod\limits_{i=1}^{l(\lambda)} P_{\lambda_i}$ and denote by $P_\lambda^* \in \left( \mathcal{C^+}\right)^*$ the basis elements dual to $P_\lambda$. We have the following.

\begin{lem}\label{lem:gluing of psi^(g,l)}
For any $g_1,g_2 \in \Z$ and $l_1, l_2 \geq 1$, we have
\begin{equation}\label{eq:gluing of two psi^(g,l)s}
\left(\sum\limits_\lambda 1^{\otimes l_1-1} \otimes  W_\lambda^* \otimes W_\lambda^* \otimes 1^{\otimes l_2-1} \right)\left(\Psi^{(-1)^{k_1}(g_1,l_1)} \otimes \Psi^{(-1)^{k_2}(g_2,l_2)} \right) = \Psi^{(-1)^{k_1+k_2}(g_1 + g_2,l_1 + l_2 - 2)}.
\end{equation}

For any $g \in \Z$ and $l \geq 2$, we have
\begin{equation}\label{eq:gluing of one psi^(g,l) in the interior}
\left(\sum\limits_{i=0}^{\infty} 1^{\otimes l-2} \otimes  W_{(i)}^* \otimes W_{(i)}^* \right) \left(\Psi^{\pm(g,l)} \right) = \Psi^{\pm(g+1,l - 2)}.
\end{equation}
\end{lem}

\bproof Let $c_\lambda \in \Q$ be the unique rational numbers so that for any power series of the form $\sum_{i=1}^\infty \frac{1}{i} A_i t^i$ with values in a $\Q$-algebra, we have that 
\begin{equation}
\exp\left(\sum\limits_{i=1}^\infty \frac{1}{i} A_i t^i \right) = \sum\limits_{i=0}^\infty \sum_{\lambda \vdash i} c_\lambda A_\lambda t^i,
\end{equation}
where $A_\lambda \coloneqq \prod_{i=1}^{l(\lambda)} A_{\lambda_i}$ for $l(\lambda) > 0$ and $A_{()} = 1$, where $()$ denotes the empty partition of $0$. Observe that $c_\lambda > 0$ for all partitions $\lambda$.

Let $\psi_1 \coloneqq \sum_{i=1}^\infty \frac{1}{i} A_i \otimes P_i$ and $\psi_2 \coloneqq \sum_{i=1}^\infty \frac{1}{i} P_i \otimes B_i$ be two power series, where $A_i$ and $B_i$ are elements of two (possibly different) commutative $\Lambda$-algebras. We will show that the partition functions $\Psi_j \coloneqq \exp \psi_j$ satisfy
\begin{equation}
\left(1 \otimes \left(\sum_\lambda W_\lambda^* \otimes W_\lambda^* \right)\otimes 1 \right)(\Psi_1 \otimes \Psi_2) = \exp \left(\sum_{i=1}^\infty \frac{1}{i} A_i \otimes B_i \right),
\end{equation}
which will imply the first claim.

First, note that 
\begin{equation}
\Psi_1 = \sum_\lambda c_\lambda A_\lambda \otimes P_\lambda
\end{equation}
and 
\begin{equation}
\Psi_2 = \sum_\lambda c_\lambda B_\lambda \otimes P_\lambda
\end{equation}
by definition of $c_\lambda$, $A_\lambda$, $B_\lambda$ and $P_\lambda$. Thus, we can write
\begin{equation}
\Psi_1 \otimes \Psi_2 = \sum_\lambda \sum_\mu c_\lambda c_\mu A_\lambda \otimes P_\lambda \otimes P_\mu \otimes B_\mu.
\end{equation}

On the other hand, define the $p(i) \times p(i)$ matrix $A$ with rational entries $A_{\lambda,\mu}$ via 
\begin{equation}
P_\lambda = \sum_\mu A_{\lambda,\mu } W_\mu,
\end{equation}
where $p(i)$ denotes the number of partitions of $i$.

Then the identity 
\begin{equation}
\sum_\lambda W_\lambda \otimes W_\lambda \overset{(\ref{eq:Psi^(0,2) in terms of W_lambda})}{=} \exp \left( \sum_{i=1}^\infty \frac{1}{i} P_i \otimes P_i \right) = \sum_{\lambda} c_\lambda P_\lambda \otimes P_\lambda = \sum_{\lambda,\mu,\nu} c_\lambda A_{\lambda,\mu} A_{\lambda,\nu} W_\mu \otimes W_\nu
\end{equation}
implies that
\begin{equation}\label{eq:w_lambda w_lambda in terms of w_mu w_nu}
\sum_\lambda c_\lambda A_{\lambda,\mu} A_{\lambda,\nu} = \delta_{\mu,\nu}\coloneqq \begin{cases}
1 \quad \text{if } \mu = \nu\\
0 \quad \text{if } \mu \neq \nu,
\end{cases}
\end{equation}
or in other words,
\begin{equation}
\widetilde{A}^T \widetilde{A} = \mathbb{I}
\end{equation}
is the identity matrix, where the matrix $\widetilde{A}$ has entries $\widetilde{A}_{\lambda,\mu} \coloneqq c_\lambda^{\frac{1}{2}} A_{\lambda,\mu}$. Therefore, we also have that $\widetilde{A} \widetilde{A}^T = \mathbb{I}$.

Using this, we compute
\begin{equation}
\begin{split}
\left( \sum_\lambda W_\lambda^* \otimes W_\lambda^* \right)(P_\mu \otimes P_\nu) &= \sum_{\rho,\eta} A_{\mu,\rho} A_{\nu,\eta} \left( \sum_\lambda W_\lambda^* \otimes W_\lambda^* \right)(W_\rho \otimes W_\eta)\\
&= \sum_\lambda A_{\mu,\lambda} A_{\nu,\lambda}\\
&= c_\mu^{-\frac{1}{2}} c_\nu^{-\frac{1}{2}} \left( \widetilde{A} \widetilde{A}^T \right)_{\mu,\nu}\\
&= c_\mu^{-1} \delta_{\mu,\nu}.
\end{split}
\end{equation}
It follows that 
\begin{equation}
\sum_\lambda W_\lambda^* \otimes W_\lambda^* = \sum_\lambda c_\lambda^{-1} P_\lambda^* \otimes P_\lambda^*.
\end{equation}

Hence,
\begin{equation}
\begin{split}
\left(1 \otimes \left(\sum_\lambda W_\lambda^* \otimes W_\lambda^* \right)\otimes 1 \right)(\Psi_1 \otimes \Psi_2) &= \sum_\lambda \sum_\mu c_\lambda c_\mu c_\lambda^{-1} \delta_{\lambda,\mu} A_\lambda \otimes B_\mu\\
&= \sum_\lambda c_\lambda A_\lambda \otimes B_\lambda\\
&= \exp \left(\sum_{i=1}^\infty \frac{1}{i} A_i \otimes B_i \right),
\end{split}
\end{equation}
as claimed.

Next, let $\psi \coloneqq \sum_{i=1}^\infty \frac{1}{i} A_i \otimes P_i \otimes P_i$ be a power series where as before the $A_i$ are elements of some commutative $\Lambda$-algebra. We will show that the partition function $\Psi \coloneqq \exp \psi$ satisfies 
\begin{equation}
\left( 1 \otimes \left(\sum_{i = 0}^\infty W_{(i)}^* \otimes W_{(i)}^* \right) \right)(\Psi) = \exp \left(\sum_{i=1}^\infty \frac{1}{i} A_i \right),
\end{equation}
which will finish the proof of Lemma \ref{lem:gluing of psi^(g,l)}. 

For this, observe that 
\begin{equation}
\Psi = \sum_\lambda c_\lambda A_\lambda \otimes P_\lambda \otimes P_\lambda.
\end{equation}
We compute that for $\lambda \vdash i$,
\begin{equation}
W_{(i)}^* \otimes W_{(i)}^* (P_\lambda \otimes P_\lambda) = A_{\lambda, (i)} A_{\lambda, (i)} = 1,
\end{equation}
where we used that $A_{\lambda,(i)} = 1$, which can be seen for instance by computing $W_{(i)}^*(P_\lambda) = \mathcal{U}(P_\lambda)|_{a=s} = 1$, where the first equality follows from the fact that $\mathcal{U}(W_\lambda)|_{a=s} = \delta_{\lambda,(i)}$ for any non-empty partition $\lambda \vdash i$.

Consequently, it follows that 
\begin{equation}
\left( 1 \otimes \left(\sum_{i = 0}^\infty W_{(i)}^* \otimes W_{(i)}^* \right) \right)(\Psi) = \sum_\lambda c_\lambda A_\lambda = \exp \left(\sum_{i=1}^\infty \frac{1}{i} A_i \right),
\end{equation}
as claimed.
\eproof\\

\section{Multicover disk crossings}\label{sec:multicover disk crossings}

In this section, we state and prove the skein valued multi-cover disk crossing formula of Theorem \ref{thm:multicover disk crossing}.

Let $\mathcal{D} \coloneqq \mathcal{S}(T^2 \setminus B)$ be the skein of a torus with a disk removed. We view $T^2 \setminus B$ as $[0,1] \times \left[\frac{1}{4},\frac{3}{4}\right] \cup \left[\frac{1}{4},\frac{3}{4}\right] \times [0,1] \subseteq [0,1] \times [0,1]$, where we identify $\{0\} \times \left[\frac{1}{4},\frac{3}{4}\right]$ with $\{1\} \times \left[\frac{1}{4},\frac{3}{4}\right]$ and $\left[\frac{1}{4},\frac{3}{4}\right] \times \{0\}$ with $\left[\frac{1}{4},\frac{3}{4}\right] \times \{1\}$. 

We can naturally embed the rectangle $\left[\frac{1}{4},\frac{3}{4}\right] \times \left[\frac{3}{4},1\right]$ into $T^2 \setminus B$. This induces a map 
\begin{equation}\label{eq:induced map by rectangle in torus without a ball}
R^n_n \otimes \mathcal{S}\left([0,1] \times \left[\frac{1}{4},\frac{3}{4}\right] \cup \left[\frac{1}{4},\frac{3}{4}\right] \times \left[0,\frac{1}{4}\right]\right) \to \mathcal{D},
\end{equation}
where $[0,1] \times \left[\frac{1}{4},\frac{3}{4}\right] \cup \left[\frac{1}{4},\frac{3}{4}\right] \times \left[0,\frac{1}{4}\right]$ has $n$ in-going boundary points along $\left[\frac{1}{4},\frac{3}{4}\right] \times \{0\}$ and $n$ out-going boundary points along $\left[\frac{1}{4},\frac{3}{4}\right] \times \{\frac{3}{4}\}$, matching the boundary points in the definition of $R^n_n$. Denote the image of this map $\mathcal{D}_n, n \in \N_0,$ and let $\mathcal{D}^+ \coloneqq \bigoplus_{n \in \N_0} \mathcal{D}_n$. This decomposition induces an $\N_0$-grading on $\mathcal{D}^+$.

\begin{thm}\label{thm:Psi^(0,1) crossing formula}
	We have that
	\begin{equation}\label{eq:disk crossing}
	\begin{tikzpicture}
	\tikzstyle{smallnode}=[circle, inner sep=0mm, outer sep=0mm, minimum size=0.5mm, draw=black, fill=black];
	\begin{knot}[
	ignore endpoint intersections=true
	clip width=5,
	clip radius=8pt,
	only when rendering/.style={
	}
	]
	\strand [thick,-{>[scale=2,length=2,width=3]}] (-1.3,0.6)
	to [out=up, in=down] (-1.3,3.6);
	\strand [thick,-{>[scale=2,length=2,width=3]}] (-2.8,2.1)
	to [out=right, in=left] (0.2,2.1);
	
	\node at (0.9,2.1) {$=$};
	\strand [thick,-{>[scale=2,length=2,width=3]}] (3.2,0.6)
	to [out=up, in=down] (3.2,3.6);
	\strand [thick,-{>[scale=2,length=2,width=3]}] (1.5,2.3)
	to [out=right, in=down, looseness=2] (2.8,3.6);
	\strand [thick] (2.8,0.6)
	to [out=up, in=left, looseness=2.1] (4.5,2.3);
	\strand [thick,-{>[scale=2,length=2,width=3]}] (1.5,1.9)
	to [out=right, in=left] (4.5,1.9);
	
	\strand (-2.8,2.7)
	to (-2.3,2.7)
	to (-2.3,3.6);
	\strand (0.2,2.7)
	to (-0.3,2.7)
	to (-0.3,3.6);
	\strand (-2.8,1.5)
	to (-2.3,1.5)
	to (-2.3,0.6);
	\strand (0.2,1.5)
	to (-0.3,1.5)
	to (-0.3,0.6);
	\strand (1.5,2.7)
	to (2,2.7)
	to (2,3.6);
	\strand (4.5,2.7)
	to (4,2.7)
	to (4,3.6);
	\strand (1.5,1.5)
	to (2,1.5)
	to (2,0.6);
	\strand (4.5,1.5)
	to (4,1.5)
	to (4,0.6);
	
	\flipcrossings{1}
	\end{knot}
	\end{tikzpicture}
	\end{equation}
	where on both sides both the vertical and the horizontal curve are decorated by $\Psi^{-(0,1)}(t)$, and the glued curve on the right hand side is decorated by
	\begin{equation}
	F^{-1}\left(\Psi^{(0,1)}(at)\right) = 1 + \sum_\lambda W_\lambda \prod_{\square \in \lambda} \frac{s^{-c(\square)}}{s^{h(\square)} - s^{-h(\square)}}.
	\end{equation}
	(I.e. we view (thickenings of) the curves as embedded solid tori and consider the induced maps from $\left(\mathcal{C}^\Lambda\right)^{\otimes 2}$, respectively $\left(\mathcal{C}^\Lambda\right)^{\otimes 3}$, to a completion of $\mathcal{D}^+$. Then the claim is that the images of $\Psi^{-(0,1)}(t) \otimes \Psi^{-(0,1)}(t)$ and $\Psi^{-(0,1)}(t) \otimes \Psi^{-(0,1)}(t) \otimes F^{-1}\left(\Psi^{(0,1)}(at)\right)$ under these two maps have the same value.)
\end{thm}

\bproof Observe that (\ref{eq:disk crossing}) takes values in a completion of $\mathcal{D}^+$ and it is enough to prove (\ref{eq:disk crossing}) degree by degree.

Note that the left hand side satisfies the relation
\begin{equation}\label{eq:disk crossing relative recursion lhs}
	\begin{tikzpicture}
		\tikzstyle{smallnode}=[circle, inner sep=0mm, outer sep=0mm, minimum size=0.7mm, draw=black, fill=black];
		\begin{knot}[
			ignore endpoint intersections=true
			clip width=5,
			clip radius=10pt,
			only when rendering/.style={
			}
			]
			\strand [thick,-{>[scale=2,length=2,width=3]}] (-1.3,0.6)
			to [out=up, in=down] (-1.3,3.6);
			\strand [thick,-{>[scale=2,length=2,width=3]}] (-2.8,2.1)
			to [out=right, in=left] (0.2,2.1);
			\node[smallnode,blue] at (-2.3,3) {};
			\strand [thick,blue,-{>[scale=2,length=2,width=3]}] (-2.3,3)
			to [out=right, in=left] (-0.3,3);
			\node[smallnode,blue] at (-0.3,3) {};
			
			\node at (0.9,2.1) {$=$};
			\strand [thick,-{>[scale=2,length=2,width=3]}] (3,0.6)
			to [out=up, in=down] (3,3.6);
			\strand [thick,-{>[scale=2,length=2,width=3]}] (1.5,2.1)
			to [out=right, in=left] (4.5,2.1);
			\node[smallnode,blue] at (2,3) {};
			\strand [thick,blue,-{>[scale=2,length=2,width=3]}] (2,3)
			to [out=right, in=left] (4,3);
			\node[smallnode,blue] at (4,3) {};
			
			\node at (5.2,2.1) {$-$};
			\strand [thick,-{>[scale=2,length=2,width=3]}] (7.3,0.6)
			to [out=up, in=down] (7.3,3.6);
			\strand [thick] (5.8,2.1)
			to [out=right, in=left] (7.1,2.1);
			\strand [thick,-{>[scale=2,length=2,width=3]}] (7.1,2.1)
			to [out=right, in=left] (8.8,2.1);
			\node[smallnode,blue] at (6.3,3) {};
			\strand [thick,blue] (6.3,3)
			to [out=right, in=down, looseness=1.8] (6.9,3.6);
			\strand [thick,blue,-{>[scale=2,length=2,width=3]}] (6.9,0.6)
			to [out=up, in=down] (6.9,1.6)
			to [out=up, in=left, looseness=2] (8.3,3);
			\node[smallnode,blue] at (8.3,3) {};
			
			\strand (-2.8,2.7)
			to (-2.3,2.7)
			to (-2.3,3.6);
			\strand (0.2,2.7)
			to (-0.3,2.7)
			to (-0.3,3.6);
			\strand (-2.8,1.5)
			to (-2.3,1.5)
			to (-2.3,0.6);
			\strand (0.2,1.5)
			to (-0.3,1.5)
			to (-0.3,0.6);
			\strand (1.5,2.7)
			to (2,2.7)
			to (2,3.6);
			\strand (4.5,2.7)
			to (4,2.7)
			to (4,3.6);
			\strand (1.5,1.5)
			to (2,1.5)
			to (2,0.6);
			\strand (4.5,1.5)
			to (4,1.5)
			to (4,0.6);
			\strand (5.8,2.7)
			to (6.3,2.7)
			to (6.3,3.6);
			\strand (8.8,2.7)
			to (8.3,2.7)
			to (8.3,3.6);
			\strand (5.8,1.5)
			to (6.3,1.5)
			to (6.3,0.6);
			\strand (8.8,1.5)
			to (8.3,1.5)
			to (8.3,0.6);
			\flipcrossings{1,2,3,4,7}
			\redraw{3}{(-2.3,3)}
			\redraw{9}{(7.3,2.1)}
			\redraw{7}{(7.3,3)}
		\end{knot}
	\end{tikzpicture}
\end{equation}
in the relative skein of $T^2 \setminus B$ with two boundary points by the relative recursion relation (\ref{eq:relative recursion for Psi^-(0,1)}) for $\Psi^{-(0,1)}$.

We show that the right hand side satisfies the same relation. This will prove (\ref{eq:disk crossing}) by induction on the degree since the two sides agree in degree $0$ and the operator $P-\bigcirc$ induced by the map (\ref{eq:induced map by rectangle in torus without a ball}) and the corresponding automorphism of $R^n_n$ is invertible on $\mathcal{D}_n$ for $n \geq 1$ as was shown in Section \ref{sec:R^n_n}. For this, we compute using the relative recursion relation for the vertical $\Psi^{-(0,1)}$ that

\begin{equation}\label{eq:disk crossing relative recursion rhs1}
	\begin{tikzpicture}
		\tikzstyle{smallnode}=[circle, inner sep=0mm, outer sep=0mm, minimum size=0.7mm, draw=black, fill=black];
		\begin{knot}[
			ignore endpoint intersections=true
			clip width=5,
			clip radius=8pt,
			only when rendering/.style={
			}
			]
			\strand [thick,-{>[scale=2,length=2,width=3]}] (-1.1,0.6)
			to [out=up, in=down] (-1.1,3.6);
			\strand [thick,-{>[scale=2,length=2,width=3]}] (-2.8,2.3)
			to [out=right, in=down, looseness=1.8] (-1.5,3.6);
			\strand [thick] (-1.5,0.6)
			to [out=up, in=left, looseness=2] (0.2,2.3);
			\strand [thick,-{>[scale=2,length=2,width=3]}] (-2.8,1.9)
			to [out=right, in=left] (0.2,1.9);
			\node[smallnode,blue] at (-2.3,3) {};
			\strand [thick,blue,-{>[scale=2,length=2,width=3]}] (-2.3,3)
			to [out=right, in=left] (-0.3,3);
			\node[smallnode,blue] at (-0.3,3) {};
			
			\node at (0.9,2.1) {$=$};
			\strand [thick,-{>[scale=2,length=2,width=3]}] (3.2,0.6)
			to [out=up, in=down] (3.2,3.6);
			\strand [thick,-{>[scale=2,length=2,width=3]}] (1.5,2.3)
			to [out=right, in=down, looseness=1.8] (2.8,3.6);
			\strand [thick] (2.8,0.6)
			to [out=up, in=left, looseness=2] (4.5,2.3);
			\strand [thick,-{>[scale=2,length=2,width=3]}] (1.5,1.9)
			to [out=right, in=left] (4.5,1.9);
			\node[smallnode,blue] at (2,3) {};
			\strand [thick,blue,-{>[scale=2,length=2,width=3]}] (2,3)
			to [out=right, in=left] (4,3);
			\node[smallnode,blue] at (4,3) {};
			
			\node at (5.2,2.1) {$-$};
			\strand [thick,-{>[scale=2,length=2,width=3]}] (7.5,0.6)
			to [out=up, in=down] (7.5,3.6);
			\strand [thick,-{>[scale=2,length=2,width=3]}] (5.8,2.3)
			to [out=right, in=down, looseness=1.8] (7.1,3.6);
			\strand [thick] (7.1,0.6)
			to [out=up, in=left, looseness=2] (8.8,2.3);
			\strand [thick,-{>[scale=2,length=2,width=3]}] (5.8,1.9)
			to [out=right, in=left] (8.8,1.9);
			\node[smallnode,blue] at (6.3,3) {};
			\strand [thick,blue] (6.3,3)
			to [out=right, in=left] (7.1,3)
			to [out=right, in=down, looseness=2] (7.7,3.6);
			\strand [thick,blue,-{>[scale=2,length=2,width=3]}] (7.7,0.6)
			to [out=up, in=down] (7.7,2.4)
			to [out=up, in=left, looseness=2] (8.3,3);
			\node[smallnode,blue] at (8.3,3) {};
			
			\strand (-2.8,2.7)
			to (-2.3,2.7)
			to (-2.3,3.6);
			\strand (0.2,2.7)
			to (-0.3,2.7)
			to (-0.3,3.6);
			\strand (-2.8,1.5)
			to (-2.3,1.5)
			to (-2.3,0.6);
			\strand (0.2,1.5)
			to (-0.3,1.5)
			to (-0.3,0.6);
			\strand (1.5,2.7)
			to (2,2.7)
			to (2,3.6);
			\strand (4.5,2.7)
			to (4,2.7)
			to (4,3.6);
			\strand (1.5,1.5)
			to (2,1.5)
			to (2,0.6);
			\strand (4.5,1.5)
			to (4,1.5)
			to (4,0.6);
			\strand (5.8,2.7)
			to (6.3,2.7)
			to (6.3,3.6);
			\strand (8.8,2.7)
			to (8.3,2.7)
			to (8.3,3.6);
			\strand (5.8,1.5)
			to (6.3,1.5)
			to (6.3,0.6);
			\strand (8.8,1.5)
			to (8.3,1.5)
			to (8.3,0.6);
			
			\flipcrossings{3,4,10,16,18,19}
			\redraw{15}{(7.1,3)}
			\redraw{11}{(7.5,2.3)}
			\redraw{16}{(7.7,2.3)}
			\redraw{13}{(7.2,1.9)}
		\end{knot}
	\end{tikzpicture}
\end{equation}

Now we apply the relative recursion relation
\begin{equation}
	\left[F^{-1}\left(\Psi^{(0,1)}(at)\right),e\right] = r\left(c,F^{-1}\left(\Psi^{(0,1)}(at)\right)\right),
\end{equation} 
which is obtained by applying an inverse full twist to the relative recursion relation (\ref{eq:variant relative recursion for Psi^(0,1)}) and adjusting the powers of $a$, to conclude that
\begin{equation}\label{eq:disk crossing relative recursion rhs1.1}
	\begin{tikzpicture}
		\tikzstyle{smallnode}=[circle, inner sep=0mm, outer sep=0mm, minimum size=0.7mm, draw=black, fill=black];
		\begin{knot}[
			ignore endpoint intersections=true
			clip width=5,
			clip radius=8pt,
			only when rendering/.style={
			}
			]
			\strand [thick,-{>[scale=2,length=2,width=3]}] (-1.1,0.6)
			to [out=up, in=down] (-1.1,3.6);
			\strand [thick,-{>[scale=2,length=2,width=3]}] (-2.8,2.3)
			to [out=right, in=down, looseness=1.8] (-1.5,3.6);
			\strand [thick] (-1.5,0.6)
			to [out=up, in=left, looseness=2] (0.2,2.3);
			\strand [thick,-{>[scale=2,length=2,width=3]}] (-2.8,1.9)
			to [out=right, in=left] (0.2,1.9);
			\node[smallnode,blue] at (-2.3,3) {};
			\strand [thick,blue,-{>[scale=2,length=2,width=3]}] (-2.3,3)
			to [out=right, in=left] (-0.3,3);
			\node[smallnode,blue] at (-0.3,3) {};
			
			\node at (0.9,2.1) {$=$};
			\strand [thick,-{>[scale=2,length=2,width=3]}] (3.2,0.6)
			to [out=up, in=down] (3.2,3.6);
			\strand [thick,-{>[scale=2,length=2,width=3]}] (1.5,2.3)
			to [out=right, in=down, looseness=1.8] (2.8,3.6);
			\strand [thick] (2.8,0.6)
			to [out=up, in=left, looseness=2] (4.5,2.3);
			\strand [thick,-{>[scale=2,length=2,width=3]}] (1.5,1.9)
			to [out=right, in=left] (4.5,1.9);
			\node[smallnode,blue] at (2,3) {};
			\strand [thick,blue,-{>[scale=2,length=2,width=3]}] (2,3)
			to [out=right, in=left] (4,3);
			\node[smallnode,blue] at (4,3) {};
			
			\node at (5.2,2.1) {$+$};
			\strand [thick,-{>[scale=2,length=2,width=3]}] (7.5,0.6)
			to [out=up, in=down] (7.5,3.6);
			\strand [thick,-{>[scale=2,length=2,width=3]}] (5.8,2.3)
			to [out=right, in=down, looseness=1.8] (7.1,3.6);
			\strand [thick] (7.1,0.6)
			to [out=up, in=left, looseness=2] (8.8,2.3);
			\strand [thick,-{>[scale=2,length=2,width=3]}] (5.8,1.9)
			to [out=right, in=left] (8.8,1.9);
			\node[smallnode,blue] at (6.3,3) {};
			\strand [thick,blue] (6.3,3)
			to [out=right, in=down, looseness=1.8] (6.9,3.6);
			\strand [thick,blue] (6.9,0.6)
			to [out=up, in=left, looseness=1.9] (8.8,2.5);
			\strand [thick,blue,-{>[scale=2,length=2,width=3]}] (5.8,2.5)
			to [out=right, in=left, looseness=2.4] (8.0,3)
			to [out=right, in=left] (8.3,3);
			\node[smallnode,blue] at (8.3,3) {};
			
			\strand (-2.8,2.7)
			to (-2.3,2.7)
			to (-2.3,3.6);
			\strand (0.2,2.7)
			to (-0.3,2.7)
			to (-0.3,3.6);
			\strand (-2.8,1.5)
			to (-2.3,1.5)
			to (-2.3,0.6);
			\strand (0.2,1.5)
			to (-0.3,1.5)
			to (-0.3,0.6);
			\strand (1.5,2.7)
			to (2,2.7)
			to (2,3.6);
			\strand (4.5,2.7)
			to (4,2.7)
			to (4,3.6);
			\strand (1.5,1.5)
			to (2,1.5)
			to (2,0.6);
			\strand (4.5,1.5)
			to (4,1.5)
			to (4,0.6);
			\strand (5.8,2.7)
			to (6.3,2.7)
			to (6.3,3.6);
			\strand (8.8,2.7)
			to (8.3,2.7)
			to (8.3,3.6);
			\strand (5.8,1.5)
			to (6.3,1.5)
			to (6.3,0.6);
			\strand (8.8,1.5)
			to (8.3,1.5)
			to (8.3,0.6);
			\flipcrossings{4,19}
			\redraw{13}{(7.2,1.9)}
			\redraw{16}{(7,1.9)}
		\end{knot}
	\end{tikzpicture}
\end{equation}

Similarly, we first use a regular isotopy and then the relative recursion relation for the horizontal $\Psi^{-(0,1)}$ to see that

\begin{equation}\label{eq:disk crossing relative recursion rhs1.2}
	\begin{tikzpicture}
		\tikzstyle{smallnode}=[circle, inner sep=0mm, outer sep=0mm, minimum size=0.7mm, draw=black, fill=black];
		\begin{knot}[
			ignore endpoint intersections=true
			clip width=5,
			clip radius=8pt,
			only when rendering/.style={
			}
			]
			\strand [thick,-{>[scale=2,length=2,width=3]}] (-1.1,0.6)
			to [out=up, in=down] (-1.1,3.6);
			\strand [thick,-{>[scale=2,length=2,width=3]}] (-2.8,2.3)
			to [out=right, in=down, looseness=1.8] (-1.5,3.6);
			\strand [thick] (-1.5,0.6)
			to [out=up, in=left, looseness=2] (0.2,2.3);
			\strand [thick,-{>[scale=2,length=2,width=3]}] (-2.8,1.9)
			to [out=right, in=left] (0.2,1.9);
			\node[smallnode,blue] at (-2.3,3) {};
			\strand [thick,blue] (-2.3,3)
			to [out=right, in=left] (-1.5,3)
			to [out=right, in=down, looseness=2] (-0.9,3.6);
			\strand [thick,blue,-{>[scale=2,length=2,width=3]}] (-0.9,0.6)
			to [out=up, in=down] (-0.9,2.4)
			to [out=up, in=left, looseness=2] (-0.3,3);
			\node[smallnode,blue] at (-0.3,3) {};
			
			\node at (0.9,2.1) {$=$};
			\strand [thick,-{>[scale=2,length=2,width=3]}] (3.2,0.6)
			to [out=up, in=down] (3.2,3.6);
			\strand [thick,-{>[scale=2,length=2,width=3]}] (1.5,2.3)
			to [out=right, in=down, looseness=1.8] (2.8,3.6);
			\strand [thick] (2.8,0.6)
			to [out=up, in=left, looseness=2] (4.5,2.3);
			\strand [thick,-{>[scale=2,length=2,width=3]}] (1.5,1.9)
			to [out=right, in=left] (4.5,1.9);
			\node[smallnode,blue] at (2,3) {};
			\strand [thick,blue] (2,3)
			to [out=right, in=down, looseness=1.8] (2.6,3.6);
			\strand [thick,blue,-{>[scale=2,length=2,width=3]}] (2.6,0.6)
			to [out=up, in=down] (2.6,1.6)
			to [out=up, in=left, looseness=2.18] (4,3);
			\node[smallnode,blue] at (4,3) {};
			
			\node at (5.2,2.1) {$=$};
			\strand [thick,-{>[scale=2,length=2,width=3]}] (7.5,0.6)
			to [out=up, in=down] (7.5,3.6);
			\strand [thick,-{>[scale=2,length=2,width=3]}] (5.8,2.3)
			to [out=right, in=down, looseness=1.8] (7.1,3.6);
			\strand [thick] (7.1,0.6)
			to [out=up, in=left, looseness=2] (8.8,2.3);
			\strand [thick] (5.8,1.9)
			to [out=right, in=left] (7,1.9);
			\strand [thick,-{>[scale=2,length=2,width=3]}] (7,1.9)
			to [out=right, in=left] (8.8,1.9);
			\node[smallnode,blue] at (6.3,3) {};
			\strand [thick,blue] (6.3,3)
			to [out=right, in=down, looseness=1.8] (6.9,3.6);
			\strand [thick,blue,-{>[scale=2,length=2,width=3]}] (6.9,0.6)
			to [out=up, in=down] (6.9,1.6)
			to [out=up, in=left, looseness=2.18] (8.3,3);
			\node[smallnode,blue] at (8.3,3) {};
			
			\node at (9.5,2.1) {$+$};
			\strand [thick,-{>[scale=2,length=2,width=3]}] (11.8,0.6)
			to [out=up, in=down] (11.8,3.6);
			\strand [thick,-{>[scale=2,length=2,width=3]}] (10.1,2.3)
			to [out=right, in=down, looseness=1.8] (11.4,3.6);
			\strand [thick] (11.4,0.6)
			to [out=up, in=left, looseness=2] (13.1,2.3);
			\strand [thick,-{>[scale=2,length=2,width=3]}] (10.1,1.9)
			to [out=right, in=left] (13.1,1.9);
			\node[smallnode,blue] at (10.6,3) {};
			\strand [thick,blue] (10.6,3)
			to [out=right, in=down, looseness=1.8] (11.2,3.6);
			\strand [thick,blue] (11.2,0.6)
			to [out=up, in=left, looseness=2] (12.3,1.7)
			to [out=right, in=left] (13.1,1.7);
			\strand [thick,blue,-{>[scale=2,length=2,width=3]}] (10.1,1.7)
			to [out=right, in=left] (10.55,1.7)
			to [out=right, in=down, looseness=2] (11.2,2.35)
			to [out=up, in=left,looseness=2] (11.85,3)
			to [out=right, in=left] (12.6,3);
			\node[smallnode,blue] at (12.6,3) {};
			
			\strand (-2.8,2.7)
			to (-2.3,2.7)
			to (-2.3,3.6);
			\strand (0.2,2.7)
			to (-0.3,2.7)
			to (-0.3,3.6);
			\strand (-2.8,1.5)
			to (-2.3,1.5)
			to (-2.3,0.6);
			\strand (0.2,1.5)
			to (-0.3,1.5)
			to (-0.3,0.6);
			\strand (1.5,2.7)
			to (2,2.7)
			to (2,3.6);
			\strand (4.5,2.7)
			to (4,2.7)
			to (4,3.6);
			\strand (1.5,1.5)
			to (2,1.5)
			to (2,0.6);
			\strand (4.5,1.5)
			to (4,1.5)
			to (4,0.6);
			\strand (5.8,2.7)
			to (6.3,2.7)
			to (6.3,3.6);
			\strand (8.8,2.7)
			to (8.3,2.7)
			to (8.3,3.6);
			\strand (5.8,1.5)
			to (6.3,1.5)
			to (6.3,0.6);
			\strand (8.8,1.5)
			to (8.3,1.5)
			to (8.3,0.6);
			\strand (10.1,2.7)
			to (10.6,2.7)
			to (10.6,3.6);
			\strand (13.1,2.7)
			to (12.6,2.7)
			to (12.6,3.6);
			\strand (10.1,1.5)
			to (10.6,1.5)
			to (10.6,0.6);
			\strand (13.1,1.5)
			to (12.6,1.5)
			to (12.6,0.6);
			
			\flipcrossings{4,6,7,15,32,34}
		\end{knot}
	\end{tikzpicture}
\end{equation}

Since the last term in (\ref{eq:disk crossing relative recursion rhs1.1}) and the last in (\ref{eq:disk crossing relative recursion rhs1.2}) are equal up to regular isotopy, it follows that 

\begin{equation}\label{eq:disk crossing relative recursion rhs final}
\begin{tikzpicture}
\tikzstyle{smallnode}=[circle, inner sep=0mm, outer sep=0mm, minimum size=0.7mm, draw=black, fill=black];
\begin{knot}[
ignore endpoint intersections=true
clip width=5,
clip radius=8pt,
only when rendering/.style={
}
]
\strand [thick,-{>[scale=2,length=2,width=3]}] (-1.1,0.6)
to [out=up, in=down] (-1.1,3.6);
\strand [thick,-{>[scale=2,length=2,width=3]}] (-2.8,2.3)
to [out=right, in=down, looseness=1.8] (-1.5,3.6);
\strand [thick] (-1.5,0.6)
to [out=up, in=left, looseness=2] (0.2,2.3);
\strand [thick,-{>[scale=2,length=2,width=3]}] (-2.8,1.9)
to [out=right, in=left] (0.2,1.9);
\node[smallnode,blue] at (-2.3,3) {};
\strand [thick,blue,-{>[scale=2,length=2,width=3]}] (-2.3,3)
to [out=right, in=left] (-0.3,3);
\node[smallnode,blue] at (-0.3,3) {};

\node at (0.9,2.1) {$=$};
\strand [thick,-{>[scale=2,length=2,width=3]}] (3.2,0.6)
to [out=up, in=down] (3.2,3.6);
\strand [thick,-{>[scale=2,length=2,width=3]}] (1.5,2.3)
to [out=right, in=down, looseness=1.8] (2.8,3.6);
\strand [thick] (2.8,0.6)
to [out=up, in=left, looseness=2] (4.5,2.3);
\strand [thick,-{>[scale=2,length=2,width=3]}] (1.5,1.9)
to [out=right, in=left] (4.5,1.9);
\node[smallnode,blue] at (2,3) {};
\strand [thick,blue,-{>[scale=2,length=2,width=3]}] (2,3)
to [out=right, in=left] (4,3);
\node[smallnode,blue] at (4,3) {};

\node at (5.2,2.1) {$-$};
\strand [thick,-{>[scale=2,length=2,width=3]}] (7.5,0.6)
to [out=up, in=down] (7.5,3.6);
\strand [thick,-{>[scale=2,length=2,width=3]}] (5.8,2.3)
to [out=right, in=down, looseness=1.8] (7.1,3.6);
\strand [thick] (7.1,0.6)
to [out=up, in=down] (7.1,1.6)
to [out=up, in=left, looseness=1.5] (7.8,2.3)
to [out=right, in=left] (8.8,2.3);
\strand [thick] (5.8,1.9)
to [out=right, in=left] (7,1.9);
\strand [thick,-{>[scale=2,length=2,width=3]}] (7,1.9)
to [out=right, in=left] (8.8,1.9);
\node[smallnode,blue] at (6.3,3) {};
\strand [thick,blue] (6.3,3)
to [out=right, in=down, looseness=1.8] (6.9,3.6);
\strand [thick,blue,-{>[scale=2,length=2,width=3]}] (6.9,0.6)
to [out=up, in=down] (6.9,2.1)
to [out=up, in=left, looseness=2] (7.8,3)
to [out=right, in=left] (8.3,3);
\node[smallnode,blue] at (8.3,3) {};

\strand (-2.8,2.7)
to (-2.3,2.7)
to (-2.3,3.6);
\strand (0.2,2.7)
to (-0.3,2.7)
to (-0.3,3.6);
\strand (-2.8,1.5)
to (-2.3,1.5)
to (-2.3,0.6);
\strand (0.2,1.5)
to (-0.3,1.5)
to (-0.3,0.6);
\strand (1.5,2.7)
to (2,2.7)
to (2,3.6);
\strand (4.5,2.7)
to (4,2.7)
to (4,3.6);
\strand (1.5,1.5)
to (2,1.5)
to (2,0.6);
\strand (4.5,1.5)
to (4,1.5)
to (4,0.6);
\strand (5.8,2.7)
to (6.3,2.7)
to (6.3,3.6);
\strand (8.8,2.7)
to (8.3,2.7)
to (8.3,3.6);
\strand (5.8,1.5)
to (6.3,1.5)
to (6.3,0.6);
\strand (8.8,1.5)
to (8.3,1.5)
to (8.3,0.6);

\flipcrossings{3,4}
\end{knot}
\end{tikzpicture}
\end{equation}
as was claimed. This finishes the proof. \eproof\\

\begin{rmk}
	By applying framing changes to the two handles of $T^2 \setminus B$, multiplying (\ref{eq:disk crossing}) from below or above by (inverse) disk partition functions, and changing signs and adjusting powers of $a$ one can obtain formulas for arbitrary crossings of (inverse) disk partition functions with arbitrary framing. For example, observe that $\mathcal{D}$ has a natural $\Z \times \Z$-grading via the intersection number with the cocores of the two handles of $T^2\setminus B$, where the $\Z$-grading corresponding to the first $\Z$-factor agrees with the $\N_0$-grading on $\mathcal{D}^+$ introduced above. Under the map $\mathcal{D} \to \mathcal{D}$ which acts as multiplication by $(-1)^{i+j}$ on elements of degree $(i,j)$, (\ref{eq:disk crossing}) is mapped to the same diagram, but now all curves are decorated by $F^{-1}\left(\Psi^{(0,1)}(at)\right)$, because the map $\mathcal{C}^+ \to \mathcal{C}^+$ which acts as multiplication by $(-1)^n$ on elements of degree $n$ maps $\Psi^{-(0,1)}(t)$ to $F^{-1}\left(\Psi^{(0,1)}(at)\right)$.
\end{rmk}

\bibliography{references}
\bibliographystyle{amsalpha}

\end{document}